\documentclass[10pt,twocolumn,twoside]{IEEEtran}
\usepackage{nomencl}
\makenomenclature
\usepackage{textcomp}
\usepackage{enumitem}
\usepackage{amsthm}
\usepackage{stackengine}
\newcommand\barbelow[1]{\stackunder[1.2pt]{$#1$}{\rule{.8ex}{.075ex}}}

\usepackage{graphicx}
\usepackage{subcaption}
\usepackage{epstopdf}
\usepackage{float}
\usepackage[nocomma]{optidef}
\usepackage{dsfont}
\usepackage{times}
\usepackage{bm}
\usepackage{fancyhdr,graphicx,amsmath,amssymb}
\usepackage[ruled,vlined]{algorithm2e}
\usepackage{hyperref}
\hypersetup{
  colorlinks   = true, 
  urlcolor     = blue, 
  linkcolor    = blue, 
  citecolor    = red   
}
\include{pythonlisting}
\usepackage{algorithmic}
\newcounter{algsubstate}


\def\bt#1{{\color{black}#1}}

\usepackage{cite}
\allowdisplaybreaks

\usepackage{xcolor}
\usepackage{amsmath}
\usepackage{amssymb}
\usepackage{mathtools}

\DeclareMathOperator*{\argmax}{arg\,max}
\makeatletter
\def\BState{\State\hskip-\ALG@thistlm}
\makeatother

\usepackage{stfloats}

\begin{document}



%
\title{A Safe First-Order Method for Pricing-Based Resource Allocation in Safety-Critical Networks}

\author{\IEEEauthorblockN{Berkay Turan \quad}
\and
\IEEEauthorblockN{Spencer Hutchinson \quad}
\and
\IEEEauthorblockN{Mahnoosh Alizadeh}%
\thanks{B. Turan, S. Hutchinson, and M. Alizadeh are with  Dept. of ECE, UCSB, Santa Barbara, CA, USA. This work is supported by  NSF grant \#1847096. E-mails: {\url{bturan@ucsb.edu},
\url{shutchinson@ucsb.edu},
\url{alizadeh@ucsb.edu}}}
}


%


\maketitle


\begin{abstract}
We introduce a novel algorithm for solving network utility maximization (NUM) problems that arise in resource allocation schemes over networks with known safety-critical constraints, where the constraints form an arbitrary convex and compact feasible set. Inspired by applications where customers' demand can only be affected through posted prices and real-time two-way communication with customers is not available, we require an algorithm to generate ``safe prices''. This means that at no iteration should the realized demand in response to the posted prices violate the safety constraints of the network. Thus, in contrast to existing distributed first-order methods, our algorithm, called safe pricing for NUM (SPNUM), is guaranteed to produce feasible primal iterates at all iterations. At the heart of the algorithm lie two key steps that must go hand in hand to guarantee safety and convergence: 1) applying a projected gradient method on a shrunk feasible set to get the desired demand, and 2) estimating the price response function of the users and determining the price so that the induced demand is close to the desired demand. We ensure safety by adjusting the shrinkage to account for the error between the induced demand and the desired demand. In addition, by gradually reducing the amount of shrinkage and the step size of the gradient method, we prove that the primal iterates produced by the SPNUM achieve a sublinear static regret of ${\cal O}(\log{(T)})$ after $T$ time steps.
\end{abstract}

%
\IEEEpeerreviewmaketitle
\newtheorem{proposition}{Proposition}
\newtheorem{definition}{Definition}
\newtheorem{corollary}{Corollary}
\newtheorem{theorem}{Theorem}
\newtheorem{lemma}{Lemma}
\newtheorem{Fact}{Fact}
\newtheorem{remark}{Remark}
\newtheorem{assumption}{Assumption}
\newtheorem{example}{Example}
\makeatletter
\def\blfootnote{\xdef\@thefnmark{}\@footnotetext}
\makeatother
\newcommand{\eqdef}{\vcentcolon=}
\newcommand{\beq}{\begin{equation}}
\newcommand{\eeq}{\end{equation}}
\newcommand{\interior}{\textnormal{int}}
\newcommand{\dom}{\textnormal{dom}}
\newcommand{\diag}{\textnormal{diag}}
\newcommand{\ie}{i.e., }
\renewcommand{\thefootnote}{\fnsymbol{footnote}}
 \renewcommand{\thefootnote}{\arabic{footnote}}
\section{Introduction}
Many applications falling within the scope of resource allocation over networks, e.g., power distribution systems \cite{samadi2010optimal}, congestion control in data networks \cite{kelly1998rate,li2023software,liu2023outlier}, wireless cellular networks \cite{chiang2004balancing}, and congestion control in urban traffic networks \cite{mehr2017joint}, deal with a multi-user optimization problem that falls under the general umbrella of \emph{network utility maximization} (NUM) problems. The shared goal in these problems is to \emph{safely} and \emph{efficiently} allocate the shared resources to the users, where safety refers to satisfying the constraints of the system that depend on the resource allocation of all the users, and efficiency refers to the total utility of the users for a given resource allocation.

In NUM problems, the user-specific utility functions are assumed to be private to the users and therefore a centralized solution is not possible. Accordingly, distributed optimization methods have become suitable tools thanks to the separable structure of NUM problems \cite{palomar2006tutorial,necoara2011parallel}. 
\bt{The idea is to decompose the main problem into sub-problems that can be solved by the individual users. The solutions of the sub-problems are then used to solve the main problem \cite{bertsekas1997nonlinear,bertsekas2015parallel}, and this has been advocated for use in different applications, e.g., \cite{li2011optimal,kelly1998rate}.} Among the two main types of decomposition methods, primal decomposition methods correspond to a direct allocation of the resources by a central coordinator and solve the primal problem, whereas dual decomposition methods based on the Lagrangian dual problem \cite{shor2012minimization} correspond to resource allocation via pricing and solve the dual problem \cite{palomar2006tutorial}. Due to the structure of NUM problems, the latter approach has been widely adopted in the literature \cite{palomar2006tutorial,nedic2009approximate,beck20141}. Additionally, it gives users the freedom of determining their own  demand based on pricing-type signals.

Although there is extensive literature on pricing algorithms based on dual decomposition, the majority of studies focus on linear constraints \cite{nedic2009approximate, beck20141, necoara2013rate, necoara2015linear, turan2022safe}, or on non-linear constraints with the assumption of separability and full user knowledge of these constraints \cite{simonetto2016primal, falsone2017dual, notarnicola2019constraint}. Furthermore, none of the aforementioned studies propose an iterative pricing algorithm that induces resource demand satisfying the hard constraints of the problem \textit{during} the iterative optimization process. Instead, these studies only provide bounds on the infeasibility amount of the resource demand (e.g., \cite{beck20141, necoara2015linear}). Our preliminary work in \cite{turan2022safe} is an exception, which is limited to problems with linear inequality constraints characterized by binary matrices. Thus, pricing-based solutions can only be realized after convergence to a near-feasible point for resource allocation systems with safety-critical constraints. Therefore, implementation of such solutions requires a negotiation process through a two-way communication network if the system has hard safety-critical constraints, which can be considered impractical in many applications.


\bt{The research presented in this paper is motivated by network resource allocation applications in safety-critical systems, where a real-time two-way communication channel with the users is not available. One particularly relevant example of this type of application can be seen in the context of pricing-based electricity demand response. When attempting to change the users' demand through posted prices, users determine their own electricity consumption to minimize their electricity bill, and no further control on the users' demand is feasible. As such, these prices must be set such that the realized demand does not violate the physical constraints of the electric grid \cite{vardakas2014survey}. This is necessary to ensure the safe and reliable operation of the grid, as violating these physical constraints could have serious reliability implications. As such, many previous works to determine prices either directly solicit all the users' preferences and solve for the prices centrally \cite{palomar2006tutorial,christakou2017ac}, or employ distributed optimization methods that require back-and-forth communication with the users to converge to an optimal and grid-safe price \cite{palomar2006tutorial,li2011optimal}. Both categories of methods have proven to be hard to implement in practical setups, motivating new research on solutions that do not require active customer engagement and still retain safe grid operations \cite{tucker2020constrained,lubin2015robust}. In light of this motivating example, the solution we devise to determine a pricing-based solution for NUM involves a number of key considerations:
\begin{enumerate}
    \item The users themselves determine their own resource demand in response to the prices, with the actual demand only becoming observable ex-post. 
    \item No negotiation or back-and-forth communication with the users is allowed, and no adjustment (curtailment) of demand is feasible, rendering existing works based on distributed optimization to determine prices inapplicable.
    \item It is essential that the safety-critical hard constraints of the systems must not be violated by users' resource demand at any time, even when their price response is unknown.
\end{enumerate}
Accordingly, the main challenge this paper aims to overcome is how to determine the prices for resources such that: 
\begin{enumerate}
    \item No preference solicitation or negotiation with the users is required.
    \item The induced resource demand of the users at every iteration always satisfies the constraints of the system (i.e., guaranteed primal feasibility).
    \item The induced resource demand of the users is efficient, i.e., the total utility earned by the users is maximized (measured through regret bounds).
\end{enumerate}
}

To this end, in this paper, we develop an iterative pricing algorithm to solve NUM problems with arbitrary convex and compact feasible sets, called safe pricing for NUM (SPNUM). We design our algorithm based solely on the realized demand in response to prices and communicate to the users only the prices for the resources at each iteration. Our contributions can be summarized as follows:
\begin{itemize}[wide]
    \item We introduce a novel algorithm, the SPNUM, for solving NUM problems with arbitrary convex and compact feasible sets through pricing. Our algorithm iteratively designs prices and allows users the freedom of determining their own decision variable based on prices according to their own profit maximization problem (without imposing any iterative variable update rule on the users).
    \item We characterize a principled way to choose algorithm parameters to guarantee feasible primal iterates at all iterations. Furthermore, we prove that the static regret incurred by the feasible primal iterates produced by the SPNUM, i.e., the cumulative gap between the optimal objective value and the objective function evaluated at the primal iterates, up to time $T$ is bounded by ${\cal O}(\log{(T)})$.
    \item We numerically evaluate our algorithm to support our theoretical findings and compare its performance to existing first-order distributed methods for NUM problems.
\end{itemize}

To the best of the authors' knowledge, no previous work has studied pricing algorithms for NUM problems on arbitrary convex feasible sets that are unknown to the users, even without consideration of safety. 
While primal-dual algorithms \cite{zhu2011distributed,koshal2011multiuser,sakurama2017distributed,turan2020resilient} can handle non-separable arbitrary convex feasible sets, they rely on a primal update rule users need to follow in order to converge as opposed to maximizing their own profit based on observed prices. To this end, our contributions extend beyond safety, since SPNUM solves NUM problems on arbitrary convex feasible sets by iteratively designing prices and allowing the users to determine their own resource demand
according to their own profit maximization problem.


The primal feasibility and the regret guarantees of the SPNUM result from a combination of two ingredients: 1) given prices and demand at a given instant, we apply a projected gradient method on a shrunk feasible set to get the next desired demand, and 2) we estimate the price response function of the users around the current prices and determine the next prices so that the induced demand is close to the desired demand. To ensure the algorithm behaves as a projected gradient method, the induced demand must be in the strict interior of the feasible set. The algorithm operates on a shrunk feasible set to account for the error between induced and desired demand, and gradually reduces shrinkage and step size to converge to the optimal solution.

\noindent
\textbf{Related work: }Besides dual (sub)gradient methods, a few other branches of literature study a similar problem to ours. We highlight how those lines of work do not meet our particular design criteria and what differentiates our work from them. Additional details on distributed optimization algorithms and their classifications can be found in the  surveys \cite{yang2019survey,zheng2022review}.
\begin{enumerate}[wide]
\item \emph{Primal-dual methods}: Primal-dual methods tackle multi-user optimization problems with arbitrary convex global constraints by applying a projected gradient descent/ascent on the primal/dual variables of the Lagrangian \cite{zhu2011distributed,koshal2011multiuser,sakurama2017distributed,turan2020resilient}. The dual variables are updated using the aggregate resource demand information of the users and can be used for pricing of the resources. Therefore the update rule for the dual variables meets our design goals. However, the primal variables, i.e., the resource demand of the users, are updated by applying one step of gradient descent instead of solving for the profit-maximizing optimal demand in response to prices. Accordingly, these algorithms do not resemble the selfish profit-maximizing behavior of the users we adopt in this paper.

\item \emph{Projected gradient methods}: The main goal of the projected gradient methods is to maintain feasibility by projecting the primal variables on the feasible convex set after each update step. Scholars have extensively studied the convergence properties of the projected gradient methods under different assumptions \cite{bertsekas1997nonlinear,calamai1987projected,levitin1966constrained}. On the other hand, the main challenge brought by our setup is that the primal variables are controlled solely by the users and cannot be manipulated (e.g., projected). Even though we can determine a feasible desired resource allocation by means of a projected gradient method, the prices that induce such resource demand are unknown due to the privacy of the utility functions, which brings unique challenges not addressed by the previous literature.

\item \emph{Interior point methods}: Interior point methods are commonly used to solve inequality-constrained problems by using barrier functions to convert them into a sequence of equality-constrained problems, which are then solved using Newton's method \cite{boyd2004convex}. While producing feasible iterates, the use of Newton's method requires the Hessian, which is often not available in practical applications, such as demand response without two-way communications. To address this limitation, previous works such as \cite{armand2000feasible} and \cite{wei2010distributed} have proposed feasible interior point methods that approximate the Hessian using first or second-order information exchange. However, these methods do not match the profit maximization rule we would like to preserve in this paper, which allows users to freely determine their resource consumption in response to posted prices.  Closest to our setup and design goals in this paper would be \cite{athuraliya2000optimization,necoara2009interior}, where separable optimization problems with linear constraints are considered. While \cite{athuraliya2000optimization} proposes a Newton-like dual update that approximates the Hessian using first-order information, only the asymptotic convergence of the algorithm is proven and the feasibility of primal iterates is not guaranteed. \cite{necoara2009interior} proposes an interior point method using Lagrangian dual decomposition with theoretical guarantees, but requires the exact Hessian for dual updates. 

\bt{
\item \emph{Constrained Online Convex Optimization}: The constrained online convex optimization literature (e.g., \cite{guo2022online,mannor2009online,yu2017online}) aims to minimize regret while establishing bounds on the constraint violation by employing iterative update algorithms on the primal variables.  A common method in constrained online convex optimization literature is updating the primal variables (i.e., resource demand of the users) directly using the gradient of the objective function as the feedback. On the contrary, our setup only allows us to update the prices and get the primal variables (i.e., resource demand of the users) as feedback afterward, where the resource demand is determined by the users according to their own profit-maximization problem. This introduces a novel challenge because both the regret and the constraints are evaluated on the primal variables, and we somehow need to set the prices such that 1) the induced demand is in the feasible region and 2) the regret incurred by the induced demand is minimized.
}
\end{enumerate}

\noindent
\textbf{Paper Organization:} The remainder of the paper is organized as follows. In Section~\ref{sec:problem}, we formalize the problem setup. In Section~\ref{sec:spnum}, we describe the SPNUM (Algorithm~\ref{alg:safenum}) and in Section~\ref{sec:regret}, we prove its feasibility and regret guarantees. In Section~\ref{sec:num}, we provide a numerical study demonstrating the efficacy of the SDGM.

\noindent
\textbf{Notation and Basic Definitions:} We denote the set of real numbers by ${\mathbb R}$ and the set of non-negative real numbers by ${\mathbb R}_+$. For vectors, $\| \cdot \|$ denotes the standard Euclidean norm and $\|\cdot\|_p$ denotes the $p$-norm. For matrices, $\|\cdot\|$ denotes the matrix norm. Given a positive integer $n>0$, $[n]$ denotes the set of integers $\{1,2,\dots,n\}$. For two vectors $x,y\in{\mathbb R}^d$, $\langle x,y\rangle$ denotes the inner product of $x$ and $y$. Given a vector $x=[x_1^\top,~x_2^\top,~\dots,~x_n^{\top}]^\top\in{\mathbb R}^d$, $x_i\in{\mathbb R}^{d_i}$ denotes the $i$'th block of $x$. For a matrix $A\in\mathbb{R}^{m\times n}$, $A_j$ denotes the $j$'th row of $A$, $A_{:,j}$ denotes the $j$'th column of $A$. Given a matrix ${A}\in {\mathbb R}^{m\times m}$, $\diag(A)\in {\mathbb R}^m$ is the vector of the diagonals of $A$, $\kappa(A)$ is the condition number of $A$, and $\sigma_{\min}(A)$/$\sigma_{\max}(A)$ are the minimum/maximum singular values of $A$. Given a function $f:{\cal X}\subseteq{\mathbb R}^d\rightarrow{\mathbb R}$, $\nabla f$ denotes the gradient of $f$, $\nabla^k f$ denotes the $k$'th order gradient of $f$, and $\dom f$ denotes the domain ${\cal X}$ of $f$. Given two vectors $x,y\in {\mathbb R}^m$, $x\leq y$ implies element-wise inequality. Given a set ${\cal X}\subset\mathbb{R}^{d}$, ${\cal X}^\interior$ denotes the interior of ${\cal X}$. Given a convex and compact set ${\cal X}\subset\mathbb{R}^{d}$ and a point $x\in{\mathbb R}^d$, ${\Pi}_{\cal X}(x)$ denotes the Euclidian projection of ${x}$ onto ${\cal X}$. We denote the closed and the open Euclidean ball with radius $r$ centered at origin as $\bar{\mathcal{B}}(r)$ and ${\cal B}(r)$,  respectively. $I_d$ denotes the identity matrix of size $d$, $\bm{1}_d$ denotes the vector of all 1's with dimension $d$, and $e_i$ denotes the unit vector with $1$ in $i$'th dimension and $0$ everywhere else. \bt{The nomenclature can be found in the Appendix.}

\begin{definition}\label{def:strong}
A differentiable function $f(\cdot)$ is said to be \textbf{$\boldsymbol{\mu}$-strongly concave} over the domain ${\cal X}$ if there exists $\mu>0$ such that
\begin{equation}\label{eq:strong}
    \langle \nabla f(x_2)-\nabla f(x_1),x_1-x_2 \rangle\geq \mu\|x_1-x_2\|^2
\end{equation}
holds for all $x_1,x_2\in \cal X$.
\end{definition}

\begin{definition}\label{def:smooth}
A differentiable function $f(\cdot)$ is said to be \textbf{$\boldsymbol{L}$-smooth} over the domain ${\cal X}$ if there exists $L>0$ such that
\begin{equation}
    \|\nabla f(x_1)-\nabla f(x_2)\|\leq L \|x_1-x_2\|
\end{equation}
holds for all $x_1,x_2\in \cal X$.
\end{definition}

\begin{definition}\label{def:lipschitz}
A function $f(\cdot)$ is said to be \textbf{$\boldsymbol{M}$-Lipschitz continuous} over the domain ${\cal X}$ if there exists $M>0$ such that
\begin{equation}
    \|f(x_1)- f(x_2)\|\leq M\|x_1-x_2\|
\end{equation}
holds for all $x_1,x_2\in \cal X$.
\end{definition}
\section{Problem Setup}\label{sec:problem}
We study the standard NUM problem \cite{kelly1998rate}, where the goal is to allocate resources to $n$ users subject to a set of coupling constraints such that the total utility of the users is maximized. It can be formulated as the following optimization problem:
\begin{subequations}\label{eq:main}
\begin{align}
    \label{eq:objective}\underset{x\in \dom f \subseteq \mathbb{R}^d}{\max}&~ f(x)=\sum_{i=1}^n f_i(x_i)\\
    \textnormal{s.t.}&~x\in {\cal X},\label{eq:constraint}
\end{align}
\end{subequations}
where $f_i(\cdot)$ is the concave utility function of user $i$ that depends on the $d_i$-dimensional vector of resource consumption, denoted by $x_i\in\dom f_i\subseteq {\mathbb R}^{d_i}$, and ${\cal X}\subset \mathbb{R}^d$ is the convex and compact set of feasible resource allocations. We also have $\sum_{i\in[n]}d_i=d$, $\dom f = \prod_{i\in[n]}\dom f_i$, and define $\bar{d}=\max_{i\in[n]}d_i$.

For all users $i\in[n]$, we define the set
    ${\cal X}_i=\{x_i\in{\mathbb R}^{d_i}:\exists x\in{\cal X} \textnormal{ s.t. } x_i \textnormal{ is the } i\textnormal{'th block of } x\}$
as the set of values that user $i$'s resource demand vector can take in the aggregate feasible set ${\cal X}$. Note that since ${\cal X}$ is convex and compact, ${\cal X}_i$ is convex and compact, ${\forall i\in[n]}$. Furthermore, if $x\in {\cal X}$, then $x_i\in {\cal X}_i$ and if $x\in {\cal X}^{\interior}$, then $x_i\in {\cal X}_i^\interior$ hold by definition. We make the following assumptions on the feasible set ${\cal X}$, and on the utility functions over ${\cal X}_i$, $\forall i\in[n]$.
\begin{assumption}\label{ass:feasibleset}
    The feasible set ${\cal X}$ is a subset of $\dom f$, i.e., ${\cal X}\subseteq \dom f$. The diameter of the feasible set ${\cal X}$ is bounded by $R$, i.e., $\|x-y\|\leq R$, $\forall x,y\in {\cal X}$. There exists a vector $\tilde{x}$ in the interior of ${\cal X}$ such that $\tilde{x}\in{\cal X}^{\interior}$.
\end{assumption}

\begin{assumption}\label{ass:utility}
For all $i\in[n]$, the utility function $f_i(\cdot)$ is $\mu$-strongly concave, $L$-smooth, $M$-Lipschitz continuous, and has $\beta$-smooth gradient over ${\cal X}_i$.
\end{assumption}
\begin{example}[Utility function]\label{ex:utility}
For instance, take $f_i(x_i)=f_\alpha(x_i)$ to be an $\alpha$-fair utility function (see
\cite{mo2000fair}) and let ${\cal X}_i=[\barbelow{x}_i,\bar{x}_i]$ with $\barbelow{x}_i>0$. We have that $\nabla f_i(x_i)\leq 1/{\barbelow{x}}_i^{\alpha}$, $-\alpha/{\barbelow{x}}_i^{\alpha+1}\leq \nabla^2 f_i^(x_i)\leq-\alpha/{\bar{x}}_i^{\alpha+1}$, and $\alpha(\alpha+1)/\bar{x}_i^{\alpha+2}\leq \nabla ^3f_i(x_i)\leq\alpha(\alpha+1)/\barbelow{x}_i^{\alpha+2}$, $\forall x\in{\cal X}_i$.  Therefore, $f_i(x_i)$ is $\alpha/{\bar{x}}_i^{\alpha+1}$-strongly concave, $\alpha/{\barbelow{x}}_i^{\alpha+1}$-smooth, and $1/{\barbelow{x}}_i^{\alpha}$-Lipschitz continuous, and has $\alpha(\alpha+1)/\barbelow{x}_i^{\alpha+2}$-smooth  gradient over ${\cal X}_i$.
\end{example}
Under Assumption~\ref{ass:utility}, the objective function \eqref{eq:objective} is strongly concave with coefficient $\mu$. Accordingly, the convex optimization problem \eqref{eq:main} has a unique solution denoted by $x^\star$ and an optimal objective value denoted by $f^\star$. 

Since $f_i(\cdot)$ are private to the users, \eqref{eq:main} cannot be solved centrally. Therefore, distributed optimization methods based on the dual decomposition framework have been proposed in the literature (e.g., \cite{palomar2006tutorial} for the case when ${\cal X}$ is a polytope) in order to incentivize selfish users with private utility functions to follow the optimal global solution. The common high-level idea is to divide the main problem into subproblems that can be solved by the individual users upon observing a pricing signal, and iteratively design prices $\{p^0,p^1,\dots\}$ to converge to the optimal resource allocation vector $x^\star$. In this framework, upon observing a price $p_i\in {\mathbb R}^{d_i}$, each user $i\in[n]$ determines their own decision variable according to their own profit maximization problem:
\begin{equation}\label{eq:priceresponse}
    g_i(p_i) = \underset{x_i\in\textnormal{dom} f_i}{\argmax}f_i(x_i) - \langle p_i,x_i\rangle.
\end{equation}
We call $g_i(\cdot)$ the price response function of user $i$ and let $g(p) = [g_1(p_1)^\top,~g_2(p_2)^\top,\dots,~g_n(p_n)^{\top}]$ be the concatenated vector of price responses given a price vector $p\in{\mathbb R}^d$.


In the next section, we propose an algorithm to iteratively design $p^t,~\forall t\geq 1$, that produce feasible primal solutions, i.e., $x^t\in{\cal X},~\forall t\geq 1$, where $x_i^t=g_i(p_i^t)$ is determined by user $i$ through \eqref{eq:priceresponse}. In addition, the algorithm should produce primal iterates that result in a sublinear static regret per user, which is measured by
\begin{equation}\label{eq:regret}
    R(T)=\frac{1}{n}\sum_{t=1}^T f^\star-f(x^t).
\end{equation}

It is worthwhile to highlight that even without the safety criterion, the literature on distributed optimization methods does not provide a distributed solution based on pricing to \eqref{eq:main} with any type of convergence guarantees. Existing works in the literature 1) utilize a pricing algorithm based on the dual decomposition framework but consider linear constraints \cite{nedic2009approximate,beck20141,necoara2013rate,necoara2015linear,turan2022safe} or non-linear and separable constraints known by the users \cite{simonetto2016primal,falsone2017dual,notarnicola2019constraint}, or 2) solve the Lagrangian dual problem by primal-dual methods \cite{zhu2011distributed,koshal2011multiuser,sakurama2017distributed,turan2020resilient}, which restrict the users to follow a primal update method that cannot be enforced in the setting where users only care about maximizing their own profit dictated by \eqref{eq:priceresponse}. Therefore, a pricing algorithm that induces a sequence of primal iterates converging to the optimal solution of \eqref{eq:main} with general convex and compact feasible sets ${\cal X}$ is novel in the distributed optimization literature.

Additionally, we note that the definition of regret in \eqref{eq:regret} quantifies the difference between the efficiencies of the optimal resource allocation and the proposed algorithm up to time $T$. When the primal iterates $\{x^t\}_{t\in[T]}$ are in the feasible set ${\cal X}$, users' resource demand can actually be realized through the posted prices without waiting for the convergence of the algorithm, and therefore regret is a meaningful measure. On the other hand, although the above sum is computable for many of the existing works mentioned earlier (e.g., \cite{beck20141,necoara2013rate} with linear constraints), they do not guarantee feasible primal iterates but only establish bounds on the amount of constraint violation at a given iteration $t$. Therefore, solutions are only realizable after convergence to a near-feasible point for resource allocation systems with safety-critical constraints. As such, they can be viewed as complex negotiations with users over what their potential demand would be in response to different prices in order to converge to the optimal price, which renders regret a less meaningful measure. By incorporating primal feasibility into our design goals, we aim to continually allocate resources to the users through posted prices \textit{during the iterative optimization process} and measure the overall efficiency of this process through regret.

\section{Safe Pricing Algorithm for NUM}\label{sec:spnum}
In this section, we describe the price update algorithm we propose, called Safe Pricing for NUM (SPNUM), that produces feasible primal iterates satisfying a sublinear regret. To do so, we will use some definitions and results from \cite{spencerl4dc} regarding the geometric properties of convex and compact sets. While the primary focus of \cite{spencerl4dc} centers on a linear stochastic bandit setup that bears little resemblance to the NUM setup under study, the definitions of the shrunk set outlined in the former are applicable to the present context as well.

\subsection{Geometric Properties of the Feasible Set}
The main ingredient that ensures the safety of SPNUM is that it operates on a shrunk feasible set, which is formally defined as follows:
\begin{definition}
\label{def:shrunk_set}
    For a compact set $\mathcal{X} \subset \mathbb{R}^d$ and a positive scalar $\Delta\in{\mathbb R}_+$, we define the \textbf{shrunk version} of $\mathcal{X}$ as $\mathcal{X}_{\Delta}:= \{ x \in \mathcal{X} : x + v \in \mathcal{X}, \forall v \in \bar{\cal B} (\Delta) \}$.
\end{definition}

\begin{example}(Shrunk polytope)
    Let $A\in {\mathbb 
    R}^{m\times d}$ and $\mathcal{X} = \{ x \in \mathbb{R}^d : Ax \leq c\}$ be a polytope. 
    The shrunk version of ${\cal X}$ is defined as ${\cal X}_{\Delta}= \{ x \in \mathbb{R}^d : A_j ^{\top} x \leq c_j - \Delta \| A_j \|,~j\in[m] \}$.
\end{example}

\begin{remark}
If ${\cal X}$ is convex and compact, then ${\cal X}_{\Delta}$ is also convex and compact.\footnote{We can equivalently define $\mathcal{X}_{\Delta}$ using Minkowski subtraction. The Minkowski subtraction of sets $A,B \subseteq \mathbb{R}^d$ is defined as $A \ominus B := \{a - b : a \in A, b \in B \}$, or equivalently, $A \ominus B = \bigcap_{b\in{B}} (A-b)$. Therefore, $\mathcal{X}_{\Delta} = \mathcal{X} \ominus \mathcal{B}(\Delta)$ is an intersection of convex and closed sets and hence is convex and closed \cite[Section 3.1]{schneider2014convex}. By Definition~\ref{def:shrunk_set}, ${\cal X}_{\Delta}$ is a subset of ${\cal X}$, and therefore bounded. A closed and bounded convex set is convex and compact.}
\end{remark}
Given the above definition of the shrunk version of a set, one can consider the maximum shrinkage that a set can withstand while still being nonempty.
We introduce the \emph{maximum shrinkage of a set} in the following definition.
\begin{definition}
\label{def:delt_max}
    For a compact set $\mathcal{X} \subset \mathbb{R}^d$, we define the \textbf{maximum shrinkage} of $\mathcal{X}$, as $H_{\mathcal{X}}:= \sup\{ \Delta : \mathcal{X}_{\Delta} \neq \emptyset \}$.
\end{definition}

\subsection{Description of the Algorithm}

\begin{algorithm}[tb]
    \caption{Safe Pricing for NUM}
    \begin{algorithmic}[1]
    \STATE {\bfseries Input:} $p^0$, $\Delta^t$, $\gamma^t$, $\eta^t$.
    \STATE \emph{(Initialization stage)}:
    \STATE \label{step:initialx}
    Each user $i\in[n]$ receives $p_i^{0}$ and $p_i^{-t} = p_i^0+\eta^0 e_{1+\mathrm{mod}(t,d_i)}$, $\forall t\in[d_i]$ and solves
    \begin{equation}\label{eq:initx}
			    x_i^{t}=g_i(p_i^t),~t=-d_i,-d_i+1,\dots, 0.
			\end{equation}
   \STATE \label{step:initialjacobian} For all $i\in[n]$, estimate the Jacobian of $g_i$ as:
   \begin{equation}
       \hat{\nabla}g_i^0=\left[\frac{x_i^{-d_i}-x_i^0}{\eta^0},\dots,~\frac{x_i^{-1}-x_i^0}{\eta^0}\right]
   \end{equation}
    \FOR{$t=0,1,\dots$}
    \STATE \emph{(Update stage)} \label{step:updatestage}
    \STATE\label{step:xupdate} Compute $\hat{x}^{t+1} = \Pi_{{\cal X}_{\Delta^t}}(x^t + \gamma^t p^t)$.
    \STATE \label{step:pupdate}Set ${p}_i^{t+1} = p_i^t +  [\hat{\nabla}g_i^t]^{-1}(\hat{x}_i^{t+1}-x_i^t)$, for all $i\in [n]$.
  \STATE Each user $i\in[n]$ receives $p_i^{t+1}$ and solves
			\begin{equation}
			    x_i^{t+1}=g_i(p_i^{t+1})
			\end{equation}
   \STATE \emph{(Sampling stage)} \label{step:samplingstage}
   \STATE \label{step:sampleprices}Each user $i\in[n]$ receives $p_i^{t+1,s}=p_i^{t+1}+\eta^{t+1}e_{1+\mathrm{mod}(t,d_i)}$ and solves
			\begin{equation}
			    x_i^{t+1,s}=g_i(p_i^{t+1,s})
			\end{equation}
   \STATE \label{step:jacobianupdate}For each user $i\in[n]$
   \begin{align}
       &[\hat{\nabla}g_i^t]_{:,1+\mathrm{mod}(t,d_i)}\leftarrow (x_i^{t+1,s}-x_i^{t+1})/{\eta^{t+1}}\\
       &\hat{\nabla}g_i^{t+1} = \hat{\nabla}g_i^t
   \end{align}
    \ENDFOR
    \end{algorithmic}
    \label{alg:safenum}
\end{algorithm}

The proposed method, called safe pricing for NUM (SPNUM) and outlined in Algorithm~\ref{alg:safenum}, consists of two stages at each iteration: 1) update stage (Step~\ref{step:updatestage}) and 2) sampling stage (Step~\ref{step:samplingstage}). The update stage proceeds similarly to a projected gradient method on the primal iterates while designing prices that induce realized iterates close to a \textit{desired iterate}. The sampling stage estimates the Jacobians of the price response functions of the users, which are used during the update stage.

In the update stage, the algorithm first determines a desired next iterate $\hat{x}^{t+1}$ in Step~\ref{step:xupdate}. However, because the primal variables are not directly controllable, prices that induce $x^{t+1}$ that is close to $\hat{x}^{t+1}$ have to be determined at Step~\ref{step:pupdate}. Accordingly, at the heart of the update stage lie two key steps: 
\begin{enumerate}[wide]
\item At iteration $t$, the central coordinator observes $x^t$ and determines the next desired iterate $\hat{x}^{t+1}$ by means of a projected gradient ascent step in Step~\ref{step:xupdate}. This is because if $x^t\in{\cal X}^\interior$, then $x_i\in{\cal X}_i^\interior$, which implies that $p_i^t = \nabla f_i(x^t)$ by Assumption~\ref{ass:utility} and the first order optimality condition for \eqref{eq:priceresponse}. Therefore, $p^t=\nabla f(x^t)$. In addition,  projection is performed onto a shrunk set ${\cal X}_{\Delta^t}$, where $\Delta^t$ controls the amount of shrinkage at time $t$. This is the key ingredient to ensure the safety of the algorithm because the uncertainty in the price response functions will cause the actual induced iterate $x^{t+1}$ in response to the price vector $p^{t+1}$ to deviate from the desired iterate $\hat{x}^{t+1}$. By adding this safety margin to the constraint, we can ensure safety if $\|x^{t+1}-\hat{x}^{t+1}\|\in\bar{\cal B}(\Delta^t)$. Finally, by utilizing a diminishing safety margin sequence $\{\Delta^t\}_{t\geq 0}$, we can ensure convergence to the optimal solution of \eqref{eq:main}.

\item Once the desired next iterate $\hat{x}^{t+1}$ is determined, the central coordinator has to determine $p_i^{t+1}$ that would ideally induce $\hat{x}_i^{t+1}$, $\forall i\in[n]$. However, the price response function is unknown to the central coordinator, and therefore an exact solution is not possible. Instead, the central coordinator makes a linear approximation of the price response function using the Jacobian estimate of $g_i$, $\forall i\in[n]$. In particular, the central coordinator keeps an estimate of the Jacobian denoted by $\hat{\nabla}g_i^t$ initialized in Steps~\ref{step:initialx} and \ref{step:initialjacobian} of the algorithm, which is constructed by varying the price vector along each dimension and estimating the gradient using the difference equation. This results in the following linear approximation of the price response function around $p_i^t$:
\begin{equation}\label{eq:ghat}
    \hat{g}_i(p) = x_i^t + \hat{\nabla}g_i^t(p-p_i^{t}).
\end{equation}
By setting $p={p}_i^{t+1}$, $\hat{g}_i({p}_i^{t+1})=\hat{x}^{t+1}$, and rearranging, we get the price update rule in Step~\ref{step:pupdate}. This requires that the $\hat{\nabla}g_i^t$ is an invertible matrix, which will be proven in Section~\ref{sec:regret}. 

\end{enumerate}

After determining $p^{t+1}$ and $x^{t+1}$, the algorithm proceeds to the sampling stage to update the Jacobian estimates. To achieve this, the central coordinator varies the price vector $p_i^{t+1}$ along the dimension $1+\mathrm{mod}(t,d_i)$ in Step~\ref{step:sampleprices} for user $i\in[n]$, resulting in a sampling price of ${p_i^{t+1,s}}$. The response is observed and denoted as $x_i^{t+1,s}$. The difference between $x_i^{t+1,s}$ and $x_i^{t+1}$ divided by the amount of price variation serves as an estimate of the gradient of the price response function along the $1+\mathrm{mod}(t,d_i)$'th principal axis, which becomes the $1+\mathrm{mod}(t,d_i)$'th column of the Jacobian estimate $\hat{\nabla}g_i^{t+1}$ in Step~\ref{step:jacobianupdate}.
It is worthwhile to highlight that for a user $i$, the error between $\hat{x}_i^{t+1}$ and $x_i^{t+1}$ has two sources: 1) the difference between the estimated Jacobian and the actual Jacobian, i.e., $\hat{\nabla}g_i^t-\nabla g_i(p_i^t)$, and 2) the high order terms not captured by the linear approximation, i.e., $R_1=g_i(p^t)-\nabla g_i(p^t)(p-p^t)$.

It is necessary that there exists an initial price vector $p^0$ such that the demand vectors in response to the initial sampling prices in \eqref{eq:initx} are in ${\cal X}^\interior$ so that the algorithm can proceed as described above. Since this has to hold before getting any feedback from the users, we make the following assumption:
\begin{assumption}
    \label{ass:initialprices}
There exists a known price vector $p^0$ such that $g(p^0)\in{\cal X}^{\interior}$ and for all $i\in[n]$, $x_i^{-d_i}\in{\cal X}_i^\interior$.
\end{assumption}
The above assumption guarantees that the initial demand vectors in \eqref{eq:initx} are in ${\cal X}^\interior_i,~\forall i\in[n]$ and therefore the initial Jacobian estimation is meaningful.

\begin{remark}\label{rem:initialset}
    One way to satisfy Assumption~\ref{ass:initialprices} is to choose $\eta^0$ such that ${\cal X}_{\frac{\sqrt{n}\eta^0}{\mu}}$ is non-empty and $p^0$ such that $g(p^0)\in {\cal X}_{\frac{\sqrt{n}\eta^0}{\mu}}$, which is proven in Appendix~\ref{app:initialset}. \bt{}
\end{remark}
\bt{
\begin{remark}\label{rem:historicalprices}
   For network resource allocation systems, the historical price response of the users can be used to choose a price point in history where the induced demand was in the feasible set. However, if there are additional assumptions that allows us to exploit the structure of the feasible set, we can utilize systematic methods. For instance, if the feasible set ${\cal X}$ is a polytope of the form ${\cal X}=\{x: Ax\leq c\}$, where $A_{ij}\geq 0$, then one way to find safe prices is to set the prices too high and gradually reduce them since low demand promotes safety. Indeed, this is a method we use in our preliminary work below (where $A$ was a binary matrix) to determine initial prices \cite{turan2022safe}.
\end{remark}
}

In the next section, we characterize a principled way to choose parameters $\Delta^t$, $\gamma^t$, and $\eta^t$ in order to produce feasible primal iterates. Additionally, we prove that the regret incurred by the iterates produced by Algorithm~\ref{alg:safenum} is ${\cal O}(\log(T))$ after $T$ iterations, and the last iterate converges to the optimal solution at the rate ${\cal O}(\log(T)/T)$.

\section{Feasibility and Regret Analysis}\label{sec:regret}
In order to prove the safety and the regret guarantees of our algorithm, we will need to bound the distance between a point in $x\in {\cal X}$ and its projection onto the shrunk set $\Pi_{{\cal X}_{\Delta}}(x)$. The following definition from \cite{spencerl4dc} formalizes this notion called the \emph{sharpness of a set}, which is defined as the maximum  distance from any point in a set to the projection of it onto the shrunk version of that set.
\begin{definition}
For a convex and compact set $\mathcal{X} \subset \mathbb{R}^d$, we define the \textbf{sharpness} of $\mathcal{X}$ as 
\begin{equation}
    \mathrm{Sharp}_{\mathcal{X}} (\Delta) \vcentcolon= \sup_{x \in {\cal X}}\|\Pi_{{\cal X}_\Delta}(x)-x\|,
\end{equation}
for all non-negative $\Delta$ such that $\mathcal{X}_{\Delta}$ is nonempty.
\end{definition}
The following proposition establishes a bound on the sharpness of convex and compact sets as a linear function of ${\Delta}$:
\begin{proposition}\cite[Corollary 11]{spencerl4dc}
\label{prop:sharp_conv}
    For a convex, compact set $\mathcal{X} \subset \mathbb{R}^d$ with non-empty interior, we have that $\mathrm{Sharp}_{\mathcal{X}} (\Delta) \leq  \Gamma_{\mathcal{X}} \Delta$ where $\Gamma_{\mathcal{X}}\geq 1$ is a constant that depends only on the geometry and the dimension of $\mathcal{X}$.
\end{proposition}
\begin{example}[Sharpness of a polytope \cite{spencerl4dc}]\label{ex:sharppolytope}
Let $\mathcal{X} = \{ x \in \mathbb{R}^d : Ax \leq c\}$ 
 be a polytope with a nonempty interior. 
 Define $\mathcal{I}_A$ to refer to the collection of all sets of $d$ indices such that for each $\{i_1, i_2, ..., i_d\} \in \mathcal{I}_A$ the vectors $A_{i_1}, A_{i_2}, ..., A_{i_d}$ are linearly independent.
For each $\ell \in \mathcal{I}_A$ where $\ell = \{i_1, i_2, ..., i_d\}$, we define $A^\ell = [A_{i_1}^\top\  A_{i_2}^\top\ ...\ A_{i_d}^\top]^\top$. 
 We have that $\mathrm{Sharp}_{\mathcal{X}}(\Delta) \leq \sqrt{d} K_{\mathcal{X}} \Delta$, where $K_{\mathcal{X}} := \max_{\ell \in \mathcal{I}_A} \kappa(A^\ell)$.

\end{example}
\begin{example}[Sharpness of a ball in ${\mathbb R}^d$]\label{ex:sharpball}
Let ${\cal X}=\{x\in\mathbb{R}^d:(x-x_0)^\top (x-x_0)\leq r^2\}$ be a ball in ${\mathbb R}^d$ with radius $r$ centered at $x_0$. We have that $\mathrm{Sharp}_{\mathcal{X}}(\Delta) = \Delta$.
\end{example} 
Although we do not specify a closed-form expression of $\Gamma_{\cal X}$ for a general convex and compact set $\cal X$, it relates to the sharpness of polytopes that are contained in ${\cal X}$, which have closed-form bounds as given by Example~\ref{ex:sharppolytope}. We refer the reader to \cite{spencerl4dc} (Proposition 10) for a detailed discussion.

The next lemma characterizes the regularity properties of $g_i(p_i)$ over the set of prices that induce a resource demand in ${\cal X}_i^\interior$ for a user $i\in[n]$. This property is crucial for our analysis and for the feasibility of the algorithm, as we need to show that the inverse of the matrix $\hat{\nabla}g_i^t$ for the price update rule in Step~\ref{step:pupdate} is a valid operation.

\begin{lemma}\label{lem:gprops}
Let ${\cal P}_i=\{p_i\in{\mathbb R}^{d_i}:g_i(p_i)\in{\cal X}_i^\interior\}$ be the set of prices that induce a resource demand in ${\cal X}_i^{\interior}$ for a user $i\in[n]$. Over ${\cal P}_i$, $g_i(p_i)$ is bijective, $1/\mu$-Lipschitz continuous, and $\beta/\mu^3$-smooth. Accordingly, $g_i(p_i)$ is invertible and $\nabla g_i(p_i) = [\nabla^2f_i(g_i(p_i))]^{-1}$.
\end{lemma}

The proof of Lemma~\ref{lem:gprops} can be found in Appendix~\ref{app:gprops}. Lemma~\ref{lem:gprops} establishes that the true Jacobian of the price response function for user $i$ is invertible because it corresponds to the inverse of the Hessian of the strongly concave utility function of user $i$. However, this does not imply that the estimated Jacobian $\hat{\nabla}g_i^t$ is invertible since it is constructed by finite difference gradient approximation. The next lemma states that the estimated Jacobian $\hat{\nabla}g_i^t$ is close enough to $\nabla g_i(p^t)$, which allows us to bound the minimum singular value of it and therefore guarantees invertibility with the appropriate choice of algorithm parameters.
\begin{lemma}\label{lem:nablagerror}
 Let $\gamma^t=1/(\mu(t+\tau))$, $\Delta^t=\Delta/(t+\tau)^2$, and $\eta^t=\mu\Delta^{t-1}/(4\sqrt{n})$ for some $\Delta>0$ and
    \begin{align}
        \label{eq:taucondition}
            \nonumber&\tau=\max\Big\{2,2\bar{d}-1,1+{2\mu\Delta\Gamma_{\cal X}}/({M\sqrt{n}}),\sqrt{{\Delta}/{{ H}_{\cal X}}},\\
            &\hspace{1cm}{L\beta M\sqrt{\bar{d}}\left(\mu+32L\Gamma_{\cal X}\sqrt{n}(\bar{d}-1)\right)}/({2\mu^4\Gamma_{\cal X}})\Big\}.
    \end{align}
     Suppose that at iteration $t$, $x^k\in{{\cal X}}_{\frac{\eta^k\sqrt{n}}{\mu}}^\interior$, $\forall k\in[\max\{t-\bar{d}+1,0\},t]$. Then, the following holds for all $i\in[n]$:
    \begin{align}
        \|\hat{\nabla}g_i^t-\nabla g_i(p_i^t)\|\leq e_i^t,
    \end{align}
    where
    \begin{equation}
        \hspace{-.1cm}e_i^t {=} \frac{2\beta\sqrt{d_i}}{\mu^3}\left(\eta^{t}{+}2L(d_i{-}1)(M\sqrt{n}\gamma^t{+}2\Delta^t\Gamma_{\cal X})\right)\leq \frac{1}{2L}.
    \end{equation}
    Accordingly, $\sigma_{\min}(\hat{\nabla}g_i^t)\geq \frac{1}{2L}$ and therefore $\hat{\nabla}g_i^t$ is invertible.
\end{lemma}
The proof of Lemma~\ref{lem:nablagerror} can be found in Appendix~\ref{app:nablagerror}. Lemma~\ref{lem:nablagerror} characterizes a principled way to choose the algorithm parameters with respect to a free parameter $\Delta$ in order to bound the difference between $\hat{\nabla}g_i^t$ and $\nabla g_i(p^t)$. In the following subsections, we will first characterize the choice of $\Delta$ that guarantees primal feasibility at all iterations and then prove the regret and convergence guarantees of Algorithm~\ref{alg:safenum} under this choice of parameters.
\subsection{Feasibility Analysis}
The following proposition characterizes the choice of the parameters $\Delta^t$, $\gamma^t$, and $\eta^t$ to ensure feasible primal iterates:
\begin{proposition}\label{prop:safety}
Let $\gamma^t = 1/(\mu(t+\tau))$ and $\Delta^t = \Delta/(t+\tau)^2$, $\eta^t=\mu\Delta^{t-1}/(4\sqrt{n})$, where $\tau$ is given  by \eqref{eq:taucondition} and
\begin{align}
    \Delta&= {\beta LMn^{3/2}(6L+\sqrt{d}(\mu/\sqrt{n}+32L(\bar{d}-1)))}/{\mu^5}.
\end{align}
Then for all $t\geq 0$, $\|\hat{x}^{t+1}-x^{t+1}\|<3\Delta^t/4$ and $\|x^{t+1}-x^{t+1,s}\|\leq\Delta^t/4$. Accordingly, for all $t\geq 0$, the iterates $x^t$ and $x^{t,s}$ produced by Algorithm~\ref{alg:safenum} are feasible and in the strict interior of the feasible set, i.e., $x^t\in{\cal X}_{\frac{\eta^t\sqrt{n}}{\mu}}^\interior$ and $x^{t,s}\in{\cal X}^{\interior},~\forall t\geq 1$.
\end{proposition}
\begin{figure*}[t]
     \centering
     \begin{subfigure}[t]{0.42\textwidth}
         \centering
         \includegraphics[width=\textwidth]{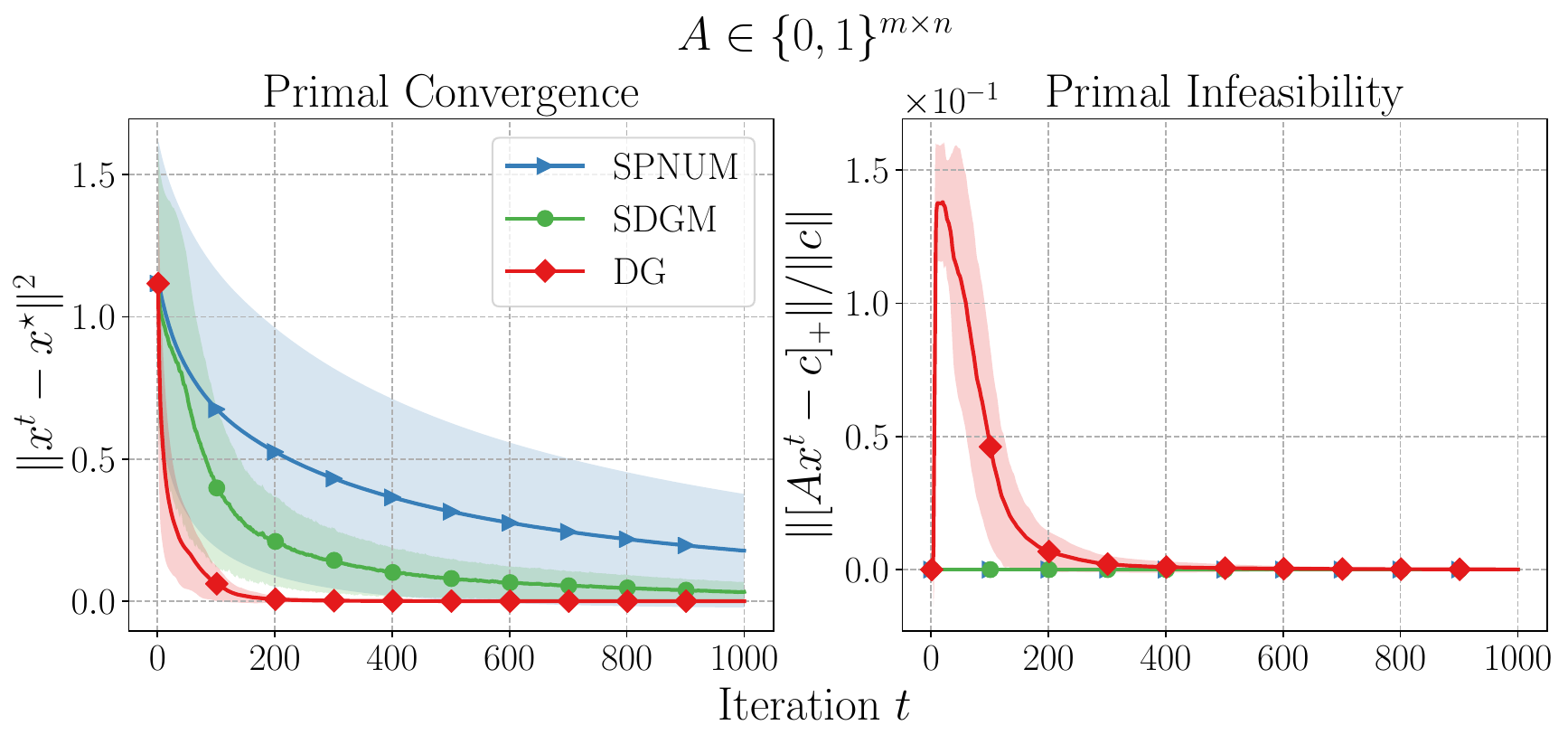}
         \vspace{-.6cm}
         \caption{}
         \label{fig:binA}
     \end{subfigure}
     \hspace{.5cm}
     \begin{subfigure}[t]{0.42\textwidth}
         \centering
         \includegraphics[width=\textwidth]{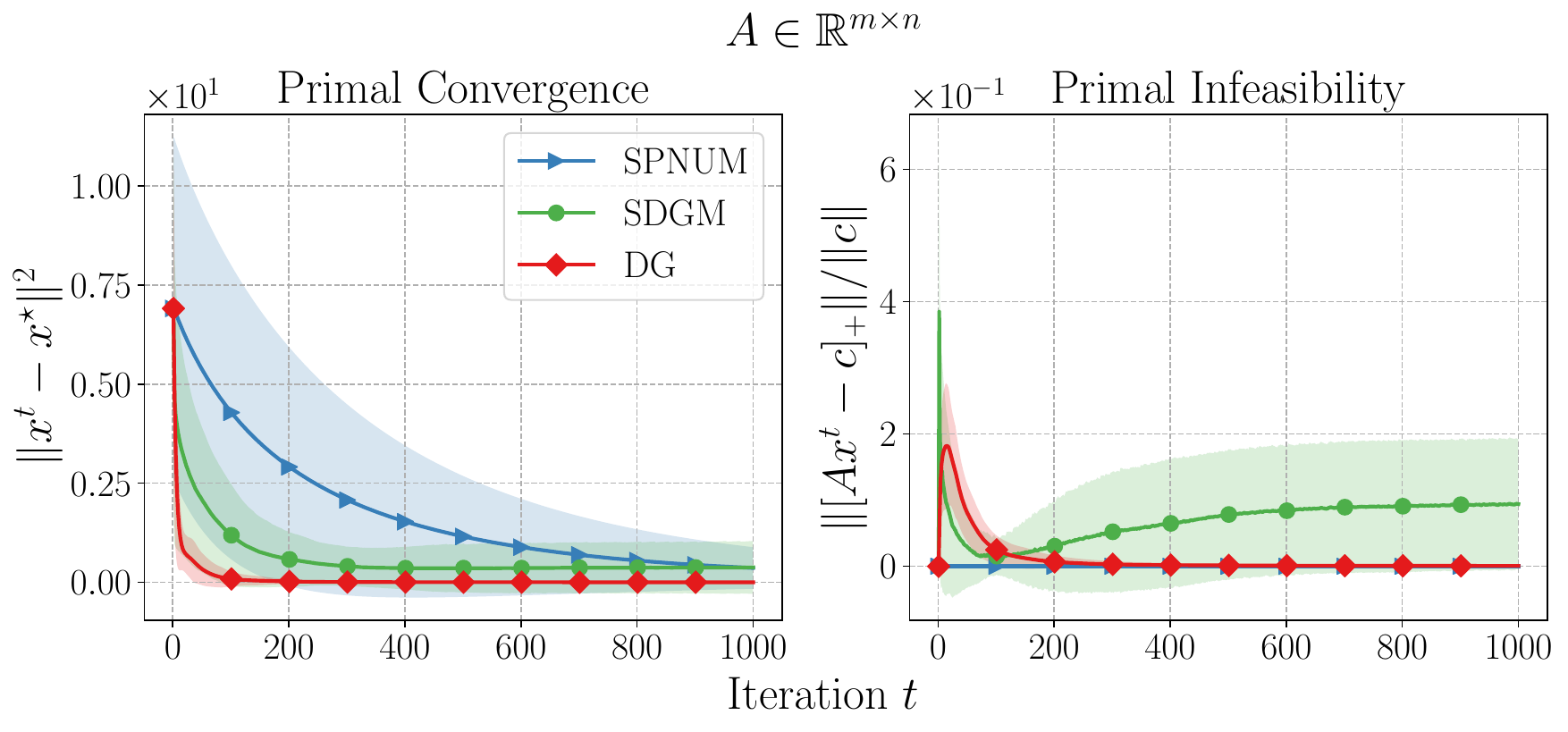}
         \vspace{-.6cm}
         \caption{}
         \label{fig:realA}
     \end{subfigure}
     \vspace{-.2cm}
        \caption{Results for the benchmarking study. In all plots, SPNUM is shown in blue, SDGM in green, and DG in red. The shaded areas correspond to one standard deviation. In (a), we plot the convergence of the primal variables measured by $\|x^t-x^\star\|^2$ (left) and the infeasibility amount measured by $\|[Ax^t-c]_+\|/\|c\|$ (right) for all three algorithms when $A\in\{0,1\}^{m\times n}$ is a binary matrix. In (b), we plot the convergence of the primal variables measured by $\|x^t-x^\star\|^2$ (left) and the infeasibility amount measured by $\|[Ax^t-c]_+\|/\|c\|$ (right) for all three algorithms when $A\in{\mathbb R}^{m\times n}$ is a real matrix.}
        \label{fig:algorithms}
        \vspace{-.5cm}
\end{figure*}
The proof of Proposition~\ref{prop:safety} can be found in Appendix~\ref{app:safety}. Given that under Proposition~\ref{prop:safety}, $x^t$ for all $t\geq 1$ are feasible and therefore implementable, the static regret \eqref{eq:regret} is a valid choice of performance metric. Next, we prove that the regret of Algorithm~\ref{alg:safenum} is ${\cal O}(\log(T))$ and the primal variables converge to the optimal solution at the rate ${\cal O}(\log(T)/T)$.
\subsection{Regret and Convergence Analysis}
As our algorithm alternates between executing one update and one sampling stage, after $T$ iterations it will have executed $T/2$ update stages and $T/2$ sampling stages. In this case, the regret per user is fairly calculated as:
\begin{equation}
    R(T)=\frac{1}{n} \sum_{t=1}^{T/2} (f(x^\star)-f(x^t)+f(x^\star)-f(x^{t,s})).
\end{equation}

The following theorem establishes an upper bound on the regret incurred by the primal iterates produced by Algorithm~\ref{alg:safenum}, and the squared distance between last iterate $x^{T/2}$ and the optimum solution $x^\star$:
\begin{theorem}\label{thm:regret}
Let $p^0$, $\gamma^t$, $\Delta^t$, and $\eta^t$ be chosen as in Proposition~\ref{prop:safety}. Then for all $t\geq 0$, the iterates produced by Algorithm~\ref{alg:safenum} are feasible. Furthermore, the regret $R(T)$ for $T\geq 2$ satisfies
\begin{equation}\label{eq:regretbound}
    R(T)\leq {\cal O}(\log(T)(1+\Delta\Gamma_{\cal X}/n)),
\end{equation}
where ${\cal O}(\cdot)$ hides other constants. In addition, the last primal iterate $x^{T/2}$ satisfies
\begin{equation}
    \|x^{T/2}-x^\star\|^2\leq {\cal O}(\log(T)/T).
\end{equation}
\end{theorem}
\textit{Proof outline:} Since the \bt{update stage of the algorithm} proceeds similarly to a projected gradient method, the proof is similar to that of a projected gradient ascent for strongly concave functions. We have an additional error term due to $\|x^{t+1}-\hat{x}^{t+1}\|$, which is ${\cal O}(\Delta^t)$. The error term impacts \bt{the regret of the update stages} as ${\cal O}(\sum_{t=1}^{T/2} \Delta^t/\gamma^t)$, which results in an additive ${\cal O}(\log(T)\Delta\Gamma_{\cal X}/n)$ term. \bt{For the regret of the sampling stages, we note that the prices for the sampling stages are set by varying the prices of the update stages by $\eta^t$. Therefore, we can upper bound the sum of the regret of all sampling stages by the regret of the update stages plus a constant additive term of $\Delta M/(4\sqrt{n})$.}

The complete proof of Theorem~\ref{thm:regret} and the explicit constants of~\eqref{eq:regretbound} can be found in Appendix~\ref{app:regret}. According to Theorem~\ref{thm:regret}, Algorithm~\ref{alg:safenum} produces feasible solutions that achieve a sublinear regret of ${\cal O}(\log(T))$. Furthermore, the primal variables induced by the prices converge to the optimal solution at the rate ${\cal O}(\log(T)/T)$.
\begin{remark}\label{rem:di1}
    When $d_i=1,~\forall i\in[n]$, $\Delta = {\cal O(}\beta n^{3/2})$ and $R(T)={\cal O}(\log(T)(1+\sqrt{n}\beta\Gamma_{\cal X}))$.
\end{remark}
In the next section, we numerically demonstrate our theoretical results about the primal variables induced by Algorithm~\ref{alg:safenum} and compare its performance to existing pricing algorithms.

\section{Numerical Studies}\label{sec:num}
In this section, we demonstrate the efficacy of SPNUM via three numerical studies: 1) a benchmarking study to compare SPNUM's convergence and feasibility performance to existing pricing methods that solve the NUM problem, 2) a toy NUM problem with a non-linear feasible set to demonstrate the success of SPNUM on non-linear feasible sets, and 3) a parameter study to demonstrate how the regret depends on the second order smoothness parameter $\beta$, sharpness parameter ${\Gamma}_{\cal X}$, and the number of users $n$.
\subsection{Benchmarking Study}

In this study, our aim is to compare the safety and convergence performance of SPNUM to the existing algorithms on feasible sets characterized by linear inequalities, i.e., ${\cal X}=\{x\in{\mathbb R}^d: Ax\leq c\}$. We compare SPNUM to DG \cite{necoara2015linear}, which can achieve a linear convergence rate, and SDGM \cite{turan2022safe}, which can provide safety when $A$ is a binary matrix.

We have implemented all algorithms on two types of $A$ matrices: 1) $A$ is a binary matrix and 2) $A$ is a real matrix. For both cases, we randomly generated a collection of 50 networks with a random number of users $n$ taking (integer) values in range $[5,20]$, and a random number of constraints $m$ taking values in the interval $[5,10]$ (generated independently). Inspired by \cite{necoara2015linear}, for all users $i\in[n]$, we let the utility function be $f_i (x_i) = -0.5 (x_i-3)^2-x_i-\theta_i \log(1+e^{x_i})$, where $\theta_i$ is sampled uniformly from $[0,1]$ for each network configuration (we shifted the quadratic term by $3$ to ensure that the optimal solution is on the boundary of the feasible set). We set $\dom f_i=[0,1]$ for all $i\in[n]$. For each network configuration, we first randomly generated a matrix $\hat{A}$ by sampling $m\times n$ Bernoulli random variables for the binary matrix case, and by sampling $m\times n$ random variables from the continuous uniform distribution in $[-1,1]$ for the real matrix case. We then let $A = [\hat{A}^\top~I_n]^\top$. For the binary case, we let $c=\bm{1}_{m+n}$, and for the real case, we let $c=[0.1\bm{1}_m^\top~\bm{1}_n^\top]^\top$.\footnote{For SPNUM, we additionally include the constraints $x\geq 0$ in ${\cal X}$ to satisfy Assumption~\ref{ass:feasibleset}. For the other algorithms, this is not needed.}

We note that ${\cal X}_i \subseteq [0,1],~\forall i\in[n]$. Within ${\cal X}_i$, using bounds on $\theta_i$ and computing the derivatives of $f_i$, we get $M=2$, $L=5/4$, $\mu=1$, $\beta=\sinh(1)/(2(1+\cosh(1))^2)\approx0.0909$. Finally, from Example~\ref{ex:sharppolytope} we have ${\Gamma}_{\cal X} \leq \sqrt{n}\kappa(A)$.

For each configuration, we initialized the dual variables and prices to induce $x_i^0=\eta^0/\mu,~\forall i\in[n]$, and ran all three algorithms for a horizon of $T=1000$. We demonstrate the results for the binary matrix case and the real matrix case in Figure~\ref{fig:binA} and Figure~\ref{fig:realA}, respectively. In Figure~\ref{fig:binA} we observe that \textbf{1)} while DG converges the fastest, it is not safe, \textbf{2)} SDGM and SPNUM converge slower but are safe, and \textbf{3)} SDGM converges faster than SPNUM because it is designed specifically for this setting. On the other hand, in Figure~\ref{fig:realA} we observe that \textbf{1)} SDGM does not provide safety and convergence when $A$ is a real matrix, as its assumptions do not hold anymore (note that the plot for $\|x^t-x^\star\|^2$ flattens for SDGM), \textbf{2)} SPNUM successfully provides safety and convergence.

\begin{figure*}[t]
    \centering
    \includegraphics[width = .8\textwidth]{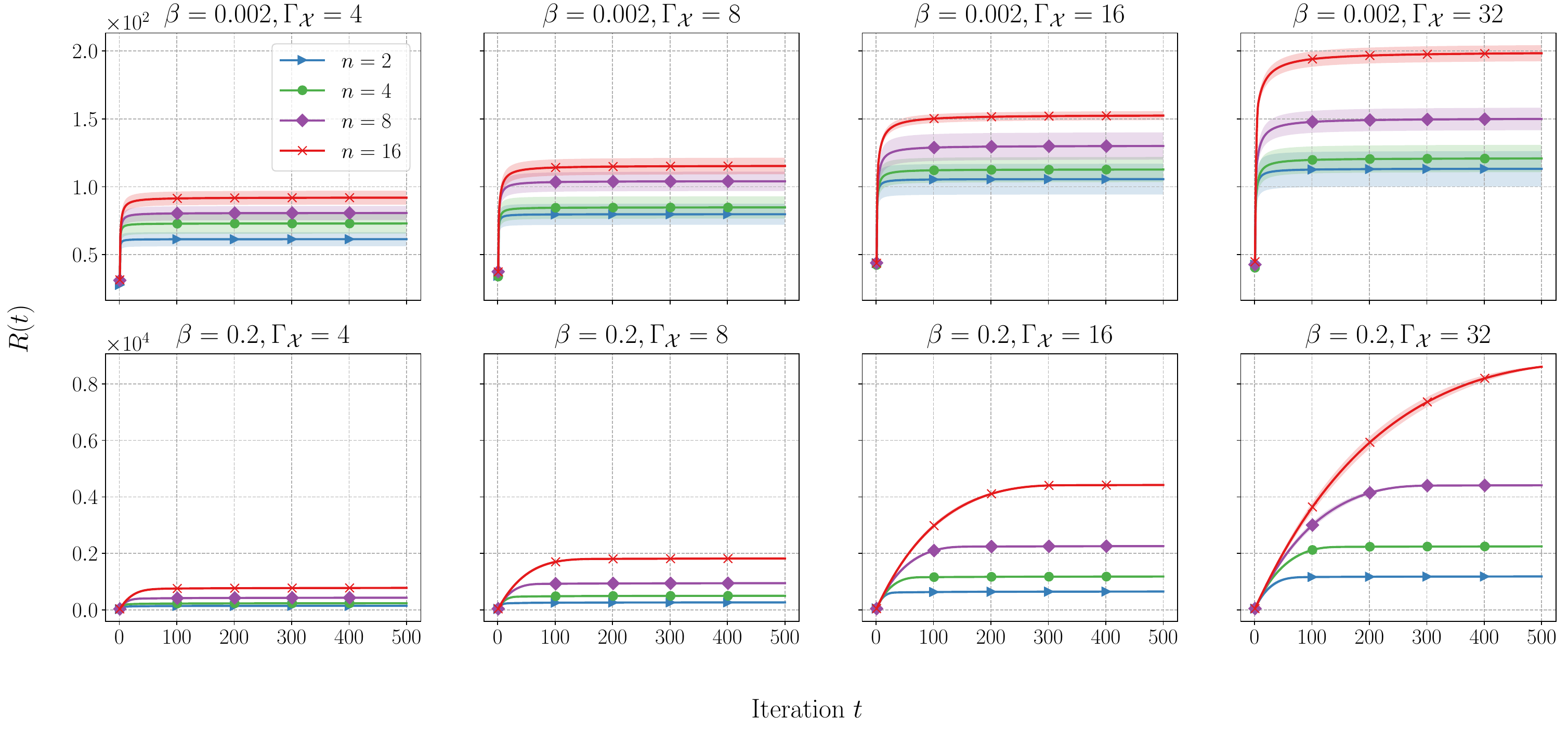}
    \vspace{-.2cm}
    \caption{Results for the numerical study on the impact of sharpness on regret. The figures on each row share the same y-axis. The shaded areas correspond to one standard deviation. The title of each subfigure denotes the $(\beta,\Gamma_{\cal X})$ configuration, and the regret incurred by different values of $n$ are plotted for each configuration. We observe that in the top row of figures, i.e., when $\beta$ is small, both $\Gamma_{\cal X}$ and $n$ have little effect on the regret (e.g., increasing ${\Gamma}_{\cal X}$ by 8 times only doubles the regret for all $n$). On the other hand, the bottom row of figures demonstrates that when $\beta$ is larger, then $\Gamma_{\cal X}$ and $n$ have a significant impact.}
    \label{fig:spnum_betagamman}
    \vspace{-.7cm}
\end{figure*}
\begin{figure}[t]
    \centering
    \includegraphics[width =.42 \textwidth]{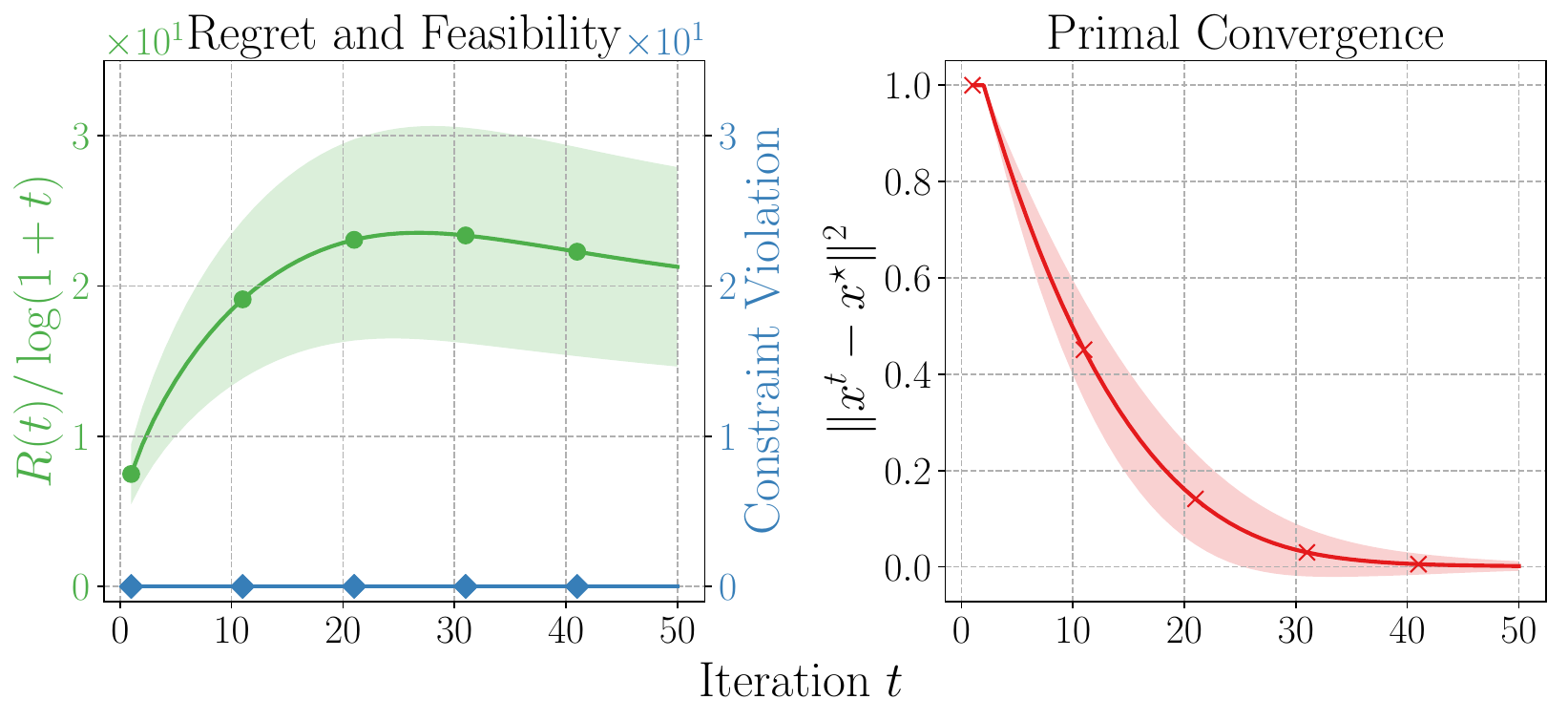}
    \vspace{-.2cm}
    \caption{Results for the numerical study on SPNUM on non-linear feasible set. In the left figure, the regret divided by ${\log(1+t)}$ is plotted in green, and constraint violation is plotted in blue, where constraint violation is $0$ if $x^t\in {\cal X}$ and $1$ otherwise. In the right figure, we plot the primal convergence measured as $\|x^t-x^\star\|^2$. Shaded areas correspond to one standard deviation.}
    \label{fig:spnum_nonlinear}
    \vspace{-.5cm}
\end{figure}
\subsection{SPNUM on Non-linear Feasible Set}

This study aims to demonstrate numerically the regret and safety guarantees of SPNUM on a problem with a feasible set characterized by non-linear inequalities. We select the feasible set ${\cal X} = \{x\in{\mathbb R}^d: \|x\|\leq 1\}$ as the unit ball in ${\mathbb R}^d$ centered at the origin. At the beginning of each run, we sample the number of users $n$ as an integer from the range $[5,20]$ uniformly at random. For all $i\in[n]$, we let the utility function be $f_i(x_i) = -0.5(x_i-y_i)^2-x_i-\theta_i\log(1+e^{x_i})$, where $\theta_i$ is sampled uniformly from $[0,1]$ and $y_i$ is sampled uniformly from $[-2,2]$ at random at the beginning of each run.

Noting that ${\cal X}_i=[-1,1]$, using bounds on $\theta_i$ and $y_i$ and computing the derivatives of $f_i$, we get $M=4+e/(1+e)$, $L=5/4$, $\mu = 1$, $\beta=\sinh(1)/(2(1+\cosh(1))^2)\approx0.0909$. Finally, from Example~\ref{ex:sharpball} we have ${\Gamma}_{\cal X} = 1$.

We initialize the prices to induce $x_i^0 = \eta^0/\mu,~\forall i\in[n]$, and ran SPNUM 100 times for a horizon of $T=50$. The results are illustrated in Figure~\ref{fig:spnum_nonlinear}. The figure shows that \textbf{1)} the regret of SPNUM grows as ${\cal O}(\log(t))$, \textbf{2)} SPNUM guarantees feasible iterates at all iterations, and \textbf{3)} the primal iterates produced by SPNUM converge to the optimal solution.

\subsection{Impact of Sharpness on Regret}

In this study, our aim is to support our theoretical results about SPNUM with numerical examples. In particular, we study the impact of sharpness parameter ${\Gamma}_{\cal X}$ and the number of users $n$ on regret through $\beta$. We set $d_i=1$, in which case $R(T) = {\cal O}(\log(T)(1+\beta\sqrt{n}\Gamma_{\cal X}))$ as stated in Remark~\ref{rem:di1}. For each user $i$, we set $f_i(x_i) = \theta_i(\cos(\omega(x-1))/\omega^2-10(x-2)^2-x\sin(\omega)/\omega)$, where $\theta_i$ is sampled uniformly from $[1,2]$. This particular choice of $f_i$ allows us to control $\beta$ while keeping the other parameters constant by simply choosing $\omega$. Using the bounds on $\theta_i$ and computing the derivatives of $f_i$, we get $M = 40$, $L=42$, $\mu=19$, and $\beta = 2\omega$.

In order to have control over the sharpness parameter ${\Gamma}_{\cal X}$, we study linear constraints of the form ${\cal X}=\{x\in \mathbb{R}^n:x\geq 0, Ax\leq c\}$, where $A_{ij}=(1-K)/(1+K(n-1))$ if $i\neq j$, and $A_{ii} = 1$. This choice of $A$ allows us to parameterize the feasible set as a function of the condition number $K$, since $\kappa(A) = K$ and ${\Gamma}_{\cal X}=\sqrt{n}\kappa(A)$. Finally, since $f_i$ is increasing over ${\cal X}_i$, the optimal solution is given by $x^\star = \bm{1}_n$.

For $n=\{2,4,8,16\}$, $\omega = \{0.001,0.1\}$, and $K=\{4/\sqrt{n},8/\sqrt{n},16/\sqrt{n},32/\sqrt{n}\}$, we randomly sampled 10 sets of $\{\theta_i\}_{i\in[n]}$, initialized $p_i^0$ so that $x_i^0 = \eta^0/\mu,~\forall i\in[n]$, and ran SPNUM for a horizon of $T=500$. Note that this corresponds to configurations of $\beta = \{0.002,0.2\}$ and $\Gamma_{\cal X}=\{4,8,16,32\}$. We plot the regret for each configuration in Figure~\ref{fig:spnum_betagamman}. The results indicate that \textbf{1)} when $\beta$ is small, $\Gamma_{\cal X}$ and $n$ have little impact on the regret, and \textbf{2)} when $\beta$ is large, regret grows with ${\Gamma}_{\cal X}$ and $n$ as the term proportional to $\sqrt{n}\beta\Gamma_{\cal X}$ dominates.

\section{Conclusion}
In this work, we introduced a novel algorithm, called the safe pricing for NUM (SPNUM), for solving resource allocation problems over networks with arbitrary convex and compact feasible sets in a distributed fashion. Our algorithm iteratively designs prices for resources and allows the users the determine their own resource demand in response to prices according to their own profit maximization problem. The prices produced by SPNUM  ensure that the induced demand satisfies the constraints of the system during the optimization process, which promotes safety. This is done by: 1) shrinking the constraint set and applying a projected gradient method to the primal variables to determine the updated desired demand, and 2) determining the prices that would induce the desired demand by estimating the price response function of the users using the historical data. By carefully controlling the amount of shrinkage to account for the error in the estimated price response, we ensure the safety of the algorithm. In addition, we have proven that the regret incurred by the SPNUM is ${\cal O}(\log(T))$, and the primal variables converge to the optimal solution at the rate of ${\cal O}(\log(T)/T)$.

\bibliographystyle{IEEEtran}
\bibliography{references}
\appendix

\subsection{Proof of Lemma~\ref{lem:nablagerror}}\label{app:nablagerror}
Firstly we note that by the choice of $\tau\geq\sqrt{\Delta/H_{\cal X}}$, we can ensure that $\Delta^t\leq H_{\cal X}$ and that ${\cal X}_{\Delta^t}$ is non-empty. Next, we show that $e_i^t\leq 1/(2L),~\forall t\geq 0$. Note that $e_i^t$ is decreasing with $t$, and therefore is maximized for $t=0$:
\begin{align}
    \hspace{-.25cm}e_i^t {\leq} e_i^0{=} {2\beta\sqrt{d_i}}\left(\eta^0{+}2L(d_i{-}1)(M\sqrt{n}\gamma^0{+}2\Delta^0\Gamma_{\cal X})\right)/{\mu^3}
\end{align}

For $\tau\geq{2\mu\Delta\Gamma_{\cal X}}/({M\sqrt{n}})$ and $d_i\leq \bar{d}$, we get:
\begin{align}
e_i^0&\leq \frac{2\beta\sqrt{\bar{d}}}{\mu^3\tau}\left(\frac{M}{8\Gamma_{\cal X}}+\frac{4L(\bar{d}-1)M\sqrt{n}}{\mu}\right)\\
&={\beta M \sqrt{\bar{d}}\left(\mu+32L\Gamma_{\cal X}\sqrt{n}(\bar{d}-1)\right)}/({4\mu^4\Gamma_{\cal X}\tau}).\label{eq:ei0bound}
\end{align}
Next, using $\tau \geq{L\beta M\sqrt{\bar{d}}\left(\mu+32L\Gamma_{\cal X}\sqrt{n}(\bar{d}-1)\right)}/({2\mu^4\Gamma_{\cal X}})$:
\begin{equation}
    e_i^t\leq e_i^0\leq 1/(2L).
\end{equation}
We will prove the lemma by induction that if $\|\hat{\nabla}g_i^k-\nabla g_i(p_i^k)\|\leq e_i^k$ holds for $k\in[\max\{0,t-d_i+1\},t-1]$, then it holds for $k = t$. 
Using Cauchy-Schwarz inequality:
\begin{align}
    \|\hat{\nabla}g_i^t{-}\nabla g_i(p_i^t)\|\leq \sqrt{d_i}\max_{j\in[d_i]}\|[\hat{\nabla}g_i^t]_{:,j}{-}[\nabla g_i(p_i^t)]_{:,j}\|.
\end{align}
For a given $j\in[d_i]$, by construction of $\hat{\nabla}g_i^t$ we have
\begin{equation}    [\hat{\nabla}g_i^t]_{:,j}=({g_i(p_i^{\ell_j}+\eta^{\ell_j} e_j)-g_i(p_i^{\ell_j})})/{\eta^{\ell_j}},
\end{equation}
for some $\ell_j\in[\max\{0,t-d_i+1\},t]$. Using the Taylor series expansion, we can rewrite the above as:
\begin{equation}
    [\hat{\nabla}g_i^t]_{:,j} = [\nabla g_i(p_i^{\ell_j})]_{:,j} + R_1/\eta^{\ell_j},
\end{equation}
where $\|R_1\|\leq \beta (\eta^{\ell_j})^2/(2\mu^3)$ follows from \cite[Lemma 1] {nesterov2006cubic} using $\beta/\mu^3$-smoothness of $g_i$. Accordingly,
\begin{align}
    \nonumber\|\hat{\nabla}g_i^t-\nabla g_i(p_i^t)\|&\leq \sqrt{d_i}\max_{j\in[d_i]}\|[\nabla g_i(p_i^{\ell_j})]_{:,j}-[\nabla g_i(p_i^t)]_{:,j}\|\\
    &~~~~+\sqrt{d_i}\beta\eta^{\ell_j}/(2\mu^3)\\
    &\hspace{-2.5cm}{\leq} \max_{{\ell_j}\in[\max\{0,t-d_i+1\},t]}\frac{\beta\sqrt{d_i}}{\mu^3}\|p_i^{\ell_j}-p_i^t\|+\frac{\sqrt{d_i}\beta\eta^{\ell_j}}{2\mu^3},\label{eq:nablagerrorcheckpt1}
\end{align}
where we used 
\begin{equation}
    \|[\nabla g_i(p_i^{\ell_j})]_{:,j}-[\nabla g_i(p_i^t)]_{:,j}\|\leq\|\nabla g_i(p_i^{\ell_j})-\nabla g_i(p_i^t)\|,
\end{equation}
for all $j\in[d_i]$, and smoothness of $g_i$. Furthermore, note that for $\tau\geq 2\bar{d}-1$, $\eta^{t-d_i+1}\leq 4\eta^t$ and therefore for $t=0$ we have
\begin{align}
    \|\hat{\nabla}g_i^0-\nabla g_i(p^0)\|&\leq 2\sqrt{d_i}\beta\eta^0/\mu^3\leq e_i^0.
\end{align}
Accordingly, the statement holds for $t=0$, which covers the base case. For $t>0$, we continue from \eqref{eq:nablagerrorcheckpt1} and bound $\|p_i^{\ell_j}-p_i^t\|$ as
\begin{align}
    \|p_i^{\ell_j}-p_i^t\|&\leq\sum_{k={\ell_j}}^{t-1}\|p_i^{k}-p_i^{k+1}\|\\
    &=\sum_{k={\ell_j}}^{t-1}\|[\hat{\nabla}g_i^k]^{-1}(\hat{x}_i^{k+1}-x_i^k)\|\\
    &\leq \sum_{k={\ell_j}}^{t-1}\|[\hat{\nabla}g_i^k]^{-1}\|\|\hat{x}_i^{k+1}-x_i^k\|.
\end{align}
The following two lemmas, whose proofs can be found in Appendices~\ref{app:sigmamin} and \ref{app:xstepamount} bound each of the terms in the above summation:
\begin{lemma}\label{lem:sigmamin}
    Suppose that $\|\hat{\nabla}g_i^t-\nabla g_i(p_i^t)\|\leq 1/(2L)$ for some $t$. Then $\sigma_{\min}(\hat{\nabla}g_i^t)\geq 1/(2L)$ and $\|[\hat{\nabla}g_i^t]^{-1}\|\leq 2L$.
\end{lemma}
\begin{lemma}\label{lem:xstepamount}
    For all $t\geq 0$, if $x^t\in{\cal X}^\interior$, then for a user $i\in[n]$ the following holds:
    \begin{equation}
        \|\hat{x}_i^{t+1}-x_i^{t}\|\leq M\sqrt{n}\gamma^t+\Delta^t\Gamma_{\cal X}.
    \end{equation}
\end{lemma}
Using Lemmas~\ref{lem:sigmamin} and \ref{lem:xstepamount}, we get
\begin{align}
   &\max_{{\ell_j}\in[\max\{0,t-d_i+1\},t]} \|p_i^{{\ell_j}}-p_i^t\|\\
   &\leq\max_{{\ell_j}\in[\max\{0,t-d_i+1\},t]}2L\sum_{k={\ell_j}}^{t-1}M\sqrt{n}\gamma^k+\Delta^k\Gamma_{\cal X}\\
   &\leq 2L(t-{\ell}_{\min})(M\sqrt{n}\gamma^{\ell_{\min}}+\Delta^{\ell_{\min}}\Gamma_{\cal X}),
\end{align}
where $\ell_{\min}=\max\{0,t-d_i+1\}$. Lastly, note that $t-\ell_{\min}\leq d_i-1$, $\gamma^{\ell_{\min}}/\gamma^t\leq 2$, and $\Delta^{\ell_{\min}}/\Delta^t\leq 4$ for $\tau\geq 2\bar{d}-1$, which gives the final result.
\subsection{Proof of Proposition~\ref{prop:safety}}\label{app:safety}
We will prove by induction that if at iteration $t$, $\forall k\in[\max\{t-\bar{d}+1,0\},t]$, $x^k\in{{\cal X}}_{\frac{\sqrt{n}\eta^k}{\mu}}^\interior$, then $x^{t+1}\in{\cal X}^{\interior}_{\frac{\sqrt{n}\eta^{t+1}}{\mu}}$ and use Assumption~\ref{ass:initialprices} that $x^{0}\in{\cal X}^{\interior}_{\frac{\sqrt{n}\eta^0}{\mu}}$. This will ensure that $x^{t+1,s}\in {\cal X}^\interior$ as well by choice of $\Delta^t$ and $\eta^t$. Therefore, we assume that $x^k\in{\cal X}^\interior_{\frac{\sqrt{n}\eta^k}{\mu}}$. Note that $\hat{x}^{t+1}\in{\cal X}^{\interior}$ by definition.

For all $i\in[n]$, we consider a modified utility function $\tilde{f}_i(x_i)$, which is equal to $f_i(x_i)$ if $x_i\in{\cal X}_i$, and an $L$-smooth, $\mu$-strongly concave extension with $\beta$-smooth gradient beyond the set ${\cal X}_i$. Accordingly, $\dom\tilde{f}_i=\mathbb{R}^{d_i}$, and $\tilde{f}_i$ is $L$-smooth and $\mu$-strongly concave over ${\mathbb R}^{d_i}$ with $\beta$-smooth gradient.

Using the modified utility function, we define the modified price response function 
\begin{equation}\label{eq:modifiedresponse}
    {\tilde g}_i(p_i)=\underset{x_i\in{\mathbb R}^{d_i}}{\argmax}\tilde{f}_i(x_i)-\langle x_i,p_i\rangle.
\end{equation}

The following Lemma, whose proof can be found in Appendix~\ref{app:smoothg}, characterizes the regularity properties of $\tilde{g}_i(p_i)$, $\forall i\in[n]$, under Assumption~\ref{ass:utility}:
\begin{lemma}\label{lem:smoothg}
    For all $i\in[n]$, let $\tilde{g}_i(p_i)$ be the modified price response function in \eqref{eq:modifiedresponse}. Then, $\tilde{g}_i(p_i)$ is bijective, $1/\mu$-Lipschitz continuous and $\beta/\mu^3$-smooth over ${\mathbb R}^{d_i}$. Furthermore, let ${\cal P}_i=\{p_i\in{\mathbb R}^{d_i}:g_i(p_i)\in{\cal X}_i^\interior\}$. The following hold true:
    \begin{enumerate}
        \item If $\tilde{g}_i(p_i)\in{\cal X}_i^\interior$, then $p_i\in{\cal P}_i$.
        \item If $p_i\in{\cal P}_i$, then $\tilde{g}_i(p_i) = g_i(p_i)$.
    \end{enumerate} 
\end{lemma}

For each user $i\in[n]$, we let $\tilde{x}_i^{t+1}=\tilde{g}_i(p_i^{t+1})$ and we rearrange the price update rule:
\begin{equation}\label{eq:xdiff}
        \tilde{x}_i^{t+1}-\hat{x}_i^{t+1} = \tilde{x}_i^{t+1}-x_i^t-\hat{\nabla}g_i^t(p^{t+1}-p^t).
    \end{equation}
    We can also write the Taylor expansion of the modified price response function $\tilde{g}_i(p)$ around $p_i^{t}$:
    \begin{equation}
        \tilde{g}_i(p_i^{t+1})-\tilde{g}_i(p_i^t)=\nabla\tilde{g}_i(p_i^t)(p_i^{t+1}-p_i^t)+R_1.
    \end{equation}
    We replace $\tilde{g}_i(p_i^{t})=g_i(p_i^t)=x_i^t$ and $\nabla\tilde{g}_i(p_i^t)=\nabla g_i(p_i^t)$ (since $p_i^t\in{\cal P}_i$) and plug the above equation into \eqref{eq:xdiff}:
     \begin{equation}
        \tilde{x}_i^{t+1}-\hat{x}_i^{t+1} = (\nabla g_i(p_i^t)-\hat{\nabla}g_i^t)(p_i^{t+1}-p_i^t) + R_1.
    \end{equation}
    To bound the norm of the above equation, we use Lemma~\ref{lem:nablagerror} to bound the norm of the first term and \cite[Lemma 1] {nesterov2006cubic} to bound the second term:
    \begin{equation}\label{eq:xdiffabs}
         \|\tilde{x}_i^{t+1}-\hat{x}_i^{t+1}\|\leq e_i^t\|p_i^{t+1}-p_i^t\|+\frac{\beta}{2\mu^3}\|p_i^{t+1}-p_i^t\|^2.
     \end{equation}
     Rearranging the price update rule and using Lemmas~\ref{lem:nablagerror} and \ref{lem:xstepamount} we can bound the norm of the price change:
      \begin{align}\label{eq:pdiffabs}
      \begin{split}
          \|p_i^{t+1}-p_i^{t}\|&\leq\|[\hat{\nabla}g_i^t]
          ^{-1}\|\|\hat{x}_i^{t+1}-x_i^t\|
      \end{split}\\
      &\leq 2L(M\sqrt{n}\gamma^t+\Delta^t\Gamma_{\cal X}).
     \end{align}
     Note that both upper bounds for $e_i^t$ and $\|p_i^{t+1}-p_i^t\|$ are decreasing with $t$. We can bound $e_i^t$ using $\tau> \frac{2\mu\Delta\Gamma_{\cal X}}{M\sqrt{n}}$ and $1\leq \Gamma_{\cal X}$ as:
     \begin{align}
         e_i^t&< {\beta M\sqrt{d_i n}(\mu/\sqrt{n}+32L(\bar{d}-1))}/({4\mu^4(t+\tau)})\\
         &={\beta M\sqrt{d_i n}\gamma^t(\mu/\sqrt{n}+32L(\bar{d}-1))}/({4\mu^3}),
     \end{align}
     and further upper bound $\|p_i^{t+1}-p_i^t\|$ as
     \begin{equation}
         \|p_i^{t+1}-p_i^t\|\leq 3LM\sqrt{n}\gamma^t.
     \end{equation}
     Plugging the above bounds and $\gamma^t$ into \eqref{eq:xdiffabs}:
     \begin{align}
     \begin{split}
          \|\tilde{x}_i^{t+1}-\hat{x}_i^{t+1}\|&< \frac{3\beta LM^2n}{4\mu^5(t+\tau)^2}\Big(6L\\
          &+{\sqrt{d_i}\left(\mu/\sqrt{n}+32L(\bar{d}-1)\right)}\Big).
     \end{split}
     \end{align}
     Next, using Cauchy-Schwarz inequality, we bound $\|\tilde{x}^{t+1}-x^{t+1}\|$ as
     \begin{align}
         \begin{split}
          \|\tilde{x}^{t+1}-\hat{x}^{t+1}\|&< \frac{3\beta LM^2n^{3/2}}{4\mu^5(t+\tau)^2}\Big(6L\\
          &+{\sqrt{d}\left(\mu/\sqrt{n}+32L(\bar{d}-1)\right)}\Big)
     \end{split}\\
     &={3\Delta^t}/{4},
     \end{align}
     where we used the definition of $\Delta^t$ and $\sum_{i\in[n]}\sqrt{d_i}\leq \sqrt{dn}$. This establishes that by definition of a shrunk set and $\Delta^t/4 = \frac{\sqrt{n}\eta^{t+1}}{\mu}$, $\tilde{x}^{t+1}\in{\cal X}^\interior_{\frac{\sqrt{n}\eta^{t+1}}{\mu}}$. Furthermore, let $\tilde{x}_i^{t+1,s}=\tilde{g_i}(p_i^{t+1,s})$. Using $1/\mu$-Lipschitz continuity of $\tilde{g_i}(p_i)$:
     \begin{equation}
         \|\tilde{x}_i^{t+1,s}-\tilde{x}_i^{t+1}\|\leq {\Delta^t}/(4\sqrt{n}),
     \end{equation}
     and $\|\tilde{x}^{t+1,s}-\tilde{x}^{t+1}\|\leq \Delta^t/4$. Accordingly, we have 
     \begin{equation}
         \|\tilde{x}^{t+1,s}-\hat{x}^{t+1}\|<\Delta^t,
     \end{equation}
     which establishes that $\tilde{x}^{t+1,s}\in{\cal X}^{\interior}$.
   Lastly, note that if $\tilde{x}^{t+1},\tilde{x}^{t+1,s}\in{\cal X}^\interior$, then for all $i\in[n]$, $\tilde{x}_i^{t+1},\tilde{x}_i^{t+1,s}\in {\cal X}_i^\interior$, or equivalently, $\tilde{g}_i(p_i^{t+1}),\tilde{g}_i(p_i^{t+1,s})\in {\cal X}_i^\interior$. Using Lemma~\ref{lem:smoothg} we have that
    $p_i^{t+1},p_i^{t+1,s}\in{\cal P}_i$, $\forall i\in[n]$. Hence, $\tilde{g}_i(p^{t+1})=g_i(p^{t+1})$ and $\tilde{x}_i^{t+1}=x_i^{t+1}$ as well as $\tilde{g}_i(p^{t+1,s})=g_i(p^{t+1,s})$ and $\tilde{x}_i^{t+1,s}=x_i^{t+1,s}$ for all $i\in [n]$, which proves the proposition.

\subsection{Proof of Theorem~\ref{thm:regret}}\label{app:regret}
We denote the regret incurred by the update stage as $R_u(T) = \sum_{t=1}^{T/2}f(x^\star)-f(x^t)$ and the regret incurred by the sampling stage as $R_s(T) = \sum_{t=1}^{T/2}f(x^\star)-f(x^{t,s})$. Let $y^{t+1}=x^t+\gamma^tp_t$. By Lemma~\ref{lem:smoothg}, we know that $p^t=\nabla f(x^t)$, $\forall t\geq 0$, since $x^t\in{\cal X}^\interior$ by Proposition~\ref{prop:safety}. For $t\geq 1$, we write using strong concavity:
\begin{align}
    &f(x^\star)-f(x^t)\leq \langle-\nabla f(x^t), x^t-x^\star\rangle-\frac{\mu}{2}\|x^t-x^\star\|^2\\
    &=\frac{1}{\gamma^t}\langle x^t-y^{t+1},x^t-x^\star\rangle-\frac{\mu}{2}\|x^t-x^\star\|^2\\
    \begin{split}\label{eq:instaregret}
            &=\frac{1}{2\gamma^t}\left(\|x^t-y^{t+1}\|^2+\|x^t-x^\star\|^2-\|y^{t+1}-x^\star\|^2\right)\\
            &\hspace{1cm}-\frac{\mu}{2}\|x^t-x^\star\|^2.
    \end{split}
\end{align}
Next, we bound the $\|y^{t+1}-x^\star\|^2$ term using Theorem~\ref{thm:schneider121} as follows:
\begin{align}
    &\|y^{t+1}-x^\star\|^2\geq \|\Pi_{{\cal X}_{\Delta^t}}(y^{t+1})-\Pi_{{\cal X}_{\Delta^t}}(x^\star)\|^2\\
    &= \|\hat{x}^{t+1}- \Pi_{{\cal X}_{\Delta^t}}(x^\star)\|^2\\
    &=\|\hat{x}^{t+1}-x^{t+1}+x^{t+1}-x^\star+x^\star-\Pi_{{\cal X}_{\Delta^t}}(x^\star)\|^2\\
           \nonumber &=\|\hat{x}^{t+1}-x^{t+1}\|^2 + \|x^{t+1}-x^\star\|^2 +\|x^\star-\Pi_{{\cal X}_{\Delta^t}}(x^\star)\|^2\hspace{-1cm}\\
          \nonumber &+2\langle\hat{x}^{t+1}{-}x^{t+1},x^{t+1}{-}x^\star\rangle {+} 2\langle x^{t+1}{-}x^\star, x^\star{-}\Pi_{{\cal X}_{\Delta^t}}(x^\star)\rangle\\
           &+2\langle x^\star-\Pi_{{\cal X}_{\Delta^t}}(x^\star),\hat{x}^{t+1}-x^{t+1}\rangle\\
    \begin{split}
           &\geq \|x^{t+1}-x^\star\|^2 - 2\|\hat{x}^{t+1}-x^{t+1}\|\|x^{t+1}-x^\star\|\\
           &-2\|x^{t+1}-x^\star\|\|x^\star-\Pi_{{\cal X}_{\Delta^t}}(x^\star)\|\\
           &-2\|x^\star-\Pi_{{\cal X}_{\Delta^t}}(x^\star)\|\|\hat{x}^{t+1}-x^{t+1}\| 
    \end{split}\\
    \begin{split}
        &\geq \|x^{t+1}-x^\star\|^2-2\Delta^tR(\Gamma_{\cal X}+3/4)-3/2(\Delta^t)^2\Gamma_{\cal X}
    \end{split}\\
    &\vcentcolon=\|x^{t+1}-x^\star\|^2-C_t,
\end{align}
where the last inequality uses $\|{x}^{t+1}-\hat{x}^{t+1}\|<3\Delta^t/4$ given by Proposition~\ref{prop:safety} and Proposition~\ref{prop:sharp_conv} to bound $\|x^\star{-}\Pi_{{\cal X}_{\Delta^t}}(x^\star)\|$. Plugging this in \eqref{eq:instaregret}:
\begin{equation}\label{eq:thm1checkpt}
\begin{split}
       f(x^\star)-f(x^t)\leq &\frac{M^2n\gamma^t}{2}-\frac{\mu}{2}\|x^t-x^\star\|^2+\frac{C^t}{2\gamma^t}\\
       &{+}\frac{1}{2\gamma^t}(\|x^t{-}x^\star\|^2-\|x^{t+1}-x^\star\|^2).
\end{split}
\end{equation}
Summing from $t=1$ to $T/2$ telescopes the $\|x^t-x^\star\|^2$ terms:
\begin{align}
\begin{split}
       & nR_u(T)\leq \frac{M^2n\log(T/2)}{2\mu}+
    \frac{\mu\tau}{2}\|x^1-x^\star\|^2\hspace{-1cm}\\
    & \hspace{1cm}+\sum_{t=2}^{T/2}\left(\frac{1}{2\gamma^t}-\frac{1}{2\gamma^{t-1}}-\frac{\mu}{2}\right)\|x^t-x^\star\|^2\\
    &\hspace{1cm}-\frac{1}{2\gamma^{T/2}}\|x^{T/2+1}-x^\star\|^2+\sum_{t=1}^{T/2}\frac{C^t}{2\gamma^t}
\end{split}\\
\label{eq:finalregret} &\leq \frac{M^2n\log(T/2)}{2\mu}+
    \frac{\mu\tau}{2}\|x^1-x^\star\|^2+\sum_{t=1}^{T/2}\frac{C^t}{2\gamma^t}.
\end{align}
Finally, note that $C^t={\cal O}(1/t^2)$ because the it consists of terms $\Delta^t$ and $(\Delta^t)^2$. Therefore, we can use the bounds ${\sum}_{t=1}^{T/2}\frac{1}{t+\tau}\leq{\sum}_{t=1}^{T/2}\frac{1}{t+2} \leq \log(T/2)$ and for $k\geq 2$, $\sum_{t=1}^{T/2}\frac{1}{(t+2)^k}\leq 1$ to show that:
\begin{equation}\label{eq:Ctsum}
\begin{split}
        \sum_{t=1}^{T/2}\frac{C^t}{2\gamma^t}\leq& \mu{\Delta R(3/4+\Gamma_{\cal X})\log(T/2)}+{3\mu\Delta^2\Gamma_{\cal X}}/4.
\end{split}
\end{equation}
Plugging \eqref{eq:Ctsum} into \eqref{eq:finalregret} and dividing by both sides by $n$, we get the regret incurred by the update stages. For the sampling stages, we note that due to the strong concavity of $f$
\begin{align}
    f(x^t){-}f(x^{t,s})&\leq\langle\nabla f(x^{t,s}),x^t{-}x^{t,s}\rangle\leq M\sqrt{n}\frac{\Delta^{t-1}}{4}.
\end{align}
Accordingly $f(x^\star)-f(x^{t,s})\leq f(x^\star)-f(x^t)+M\sqrt{n}\Delta^{t-1}/4$. Summing from $t=1$ to $T/2$, we get
\begin{align}
    nR_s(T)&{=}nR_u(T){+}\frac{M}{4}\sum_{t=1}^{T/2}\Delta^{t-1}{\leq} nR_u(T){+}\frac{\Delta M\sqrt{n}}{4},
\end{align}
which gives the final result as
\begin{equation}
    R(T)\leq 2R_u(T)+{\Delta M}/({4\sqrt{n}}).
    \end{equation}
To get the convergence result, we rearrange \eqref{eq:thm1checkpt}:
\begin{align}
\begin{split}
        \|x^{t+1}-x^\star\|^2&\leq \|x^t-x^\star\|^2(1-\mu\gamma^t) +M^2 n (\gamma^t)^2\\
        & + C^t +2\gamma^t (f(x^t)-f(x^\star))
\end{split}\\
\leq& \|x^t-x^\star\|^2(1-\mu\gamma^t) +M^2 n (\gamma^t)^2+ C^t.
\end{align}
We get an equation like the above for all $t\geq 0$. We multiply each by $(1-\mu\gamma^{t+1})$ for $t<T/2-1$ and sum them from $t=0$ to $t= T/2-1$ to get:
\begin{align}
    \begin{split}
        &\|x^{T/2}-x^\star\|^2\leq \|x^0-x^\star\|^2\prod_{t=0}^{T/2-1}(1-\mu\gamma^t)\\
        &+M^2n\sum_{t=0}^{T/2-1}(\gamma^t)^2\prod_{i=t+1}^{T/2-1}(1-\mu\gamma^i)\\
        &+\sum_{t=0}^{T/2-1}C^t\prod_{i=t+1}^{T/2-1}(1-\mu\gamma^i)\hspace{-1cm}
    \end{split}\\
    \begin{split}
        &\leq \|x^0-x^\star\|^2\frac{\tau-1}{\tau-1+T/2}+\frac{M^2n\log(T/2)}{\mu^2(T/2+\tau-1)}\\
    &+\frac{2R(3/4+\Gamma_{\cal X})\Delta\log(T/2)}{(T/2+\tau-1)}+\frac{3\Delta^2\Gamma_{\cal X}}{2(T/2+\tau-1)}.
    \end{split}
\end{align}
which completes the proof.

\subsection{Proof of Remark~\ref{rem:initialset}}\label{app:initialset}
For a user $i\in[n]$, using the modified price response function $\tilde{g}_i(p_i)$ introduced in the proof of Proposition~\ref{prop:safety}, we have that
\begin{equation}
    \|\tilde{x}_i^{-t}-x_i^0\|\leq{\eta^0}/{\mu},~\forall t\in[-d_i,-1],
\end{equation}
which implies that $\tilde{x}_i^{-t}\in{\cal X}_i^\interior$ because $x^0\in{\cal X}^\interior_{\frac{\eta^0\sqrt{n}}{\mu}}$. As such, $\tilde{x}_i^{-t} = x_i^{-t}$ and $p_i^{-t}= \nabla f_i(x_i^{-t})$.

\begin{IEEEbiographynophoto}{BERKAY TURAN}is pursuing the Ph.D. degree in Electrical and Computer Engineering at the University of California, Santa Barbara. He received the B.Sc. degree in Electrical and Electronics Engineering as well as the B.Sc. degree in  Physics degree from \ Bo\u gazi\c ci University, Istanbul, Turkey, in 2018. The overarching goal of his research is to design network control, optimization, and learning frameworks to promote efficiency and resiliency in societal-scale cyber-physical systems.
\end{IEEEbiographynophoto}
\vskip -2\baselineskip plus -1fil
\begin{IEEEbiographynophoto}{SPENCER HUTCHINSON} received the B.S. degree in electrical engineering from Colorado School of Mines in 2021. He is currently pursuing the Ph.D. degree in electrical and computer engineering from the University of California, Santa Barbara in Santa Barbara, CA, USA. His research interests include the design and analysis of optimization and learning algorithms for the control of human-cyber-physical systems.
\end{IEEEbiographynophoto}
\vskip -2\baselineskip plus -1fil
\begin{IEEEbiographynophoto}{MAHNOOSH ALIZADEH} is an associate professor of Electrical and Computer Engineering at the University of California Santa Barbara. 
She received the B.Sc. degree (’09) in Electrical Engineering from Sharif University of Technology and the M.Sc. (’13) and Ph.D. (’14) degrees in Electrical and Computer Engineering from the University of California Davis. From 2014 to 2016, she was a postdoctoral scholar at Stanford University. Her research interests are focused on designing network control, optimization, and learning frameworks to promote efficiency and resiliency in societal-scale cyber-physical systems. Dr. Alizadeh is a recipient of the NSF CAREER award.
\end{IEEEbiographynophoto}
\providetoggle{nomsort}
\settoggle{nomsort}{true} 
\newpage
\subsection{Proof of Lemma~\ref{lem:gprops}}\label{app:gprops}
By definition, $f_i(x_i)$ is strongly concave over ${\cal X}_i$, therefore the optimization problem ${\max}_{x\in{\dom f_i}} f_i(x_i)-\langle x_i,p_i\rangle$ is strongly concave and has a unique solution for $p_i\in{\cal P}_i$. \bt{Since ${\cal X}_i\subseteq \dom f_i$ by Assumption~\ref{ass:feasibleset}, the optimal solution is in the interior of $\dom f_i$.} Therefore the first-order optimality condition implies that the optimal solution $g_i(p_i)$ satisfies
\begin{equation}
    p_i = \nabla f_i(g_i(p_i)),
\end{equation}
which implies that $\nabla{f}_i$ is surjective for $p_i\in{\cal P}_i$. We also know that the gradient of a strongly concave function is injective\footnote{To see this, suppose that $x_1\neq x_2$ and therefore $\|x_1-x_2\|>0$. If $\nabla f(x_1)= \nabla f(x_2)$, \eqref{eq:strong} results in $0\geq \mu\|x_1-x_2\|^2$, which is a contradiction and $x_1=x_2$ must hold.}, therefore, $\nabla {f}_i$ is bijective and invertible and ${g}_i(p_i)=\nabla{f}_i^{-1}(p_i)$, which also proves that $g_i(p_i)$ is bijective. By the inverse function theorem, we get that:
\begin{equation}
    \nabla {g}_i(p_i) = [\nabla^2 {f}_i({g}_i(p_i))]^{-1}.
\end{equation}
Since ${f}_i$ is $L$-smooth and $\mu$-strongly concave, inverse of it's Hessian has eigenvalues in $[-1/\mu,-1/L]$, which results in
\begin{equation}
    \| \nabla {g}_i(p_i)\| = \|[\nabla^2 {f}_i({g}_i(p_i))]^{-1}\| \leq 1/\mu,
\end{equation}
proving the Lipschitz property of ${g}_i(p_i)$. To show smoothness, we let $x_i^1= {g}_i(p_i^1)$ and $x_i^2={g}_i(p_i^2)$ and write:
\begin{align}
   &\hspace*{-.2cm}\|\nabla {g}_i(p_i^1){-}\nabla {g}_i(p_i^2)\|=\|[\nabla^2 {f}_i(x_i^1)]^{-1}{-}[\nabla^2 {f}_i(x_i^2)]^{-1}\|\\
   &\hspace*{-.2cm}=\|[\nabla^2 {f}_i(x_i^1)]^{-1}(\nabla^2 {f}_i(x_i^2){-}\nabla^2 {f}_i(x_i^1))[\nabla^2 {f}_i(x_i^2)]^{-1}\|\\
   &\hspace*{-.2cm}\leq {\beta}\|x_i^1-x_i^2\|/{\mu^2}\leq {\beta}\|p_i^1-p_i^2\|/{\mu^3},
\end{align}
where the last inequality uses $1/\mu$-Lipschitz continuity of ${g}_i(p_i)$, which proves $\beta/\mu^3$-smoothness of ${g}_i(p_i)$.
\subsection{Proof of Lemma~\ref{lem:sigmamin}}\label{app:sigmamin}
Note that for $p_i^t\in{\cal P}_i$, $\nabla g_i(p^t)=[\nabla^2 f_i(g_i(p^t))]^{-1}$ is symmetric by Schwarz's theorem, since $\nabla^2 f_i(g_i(p_i))$ is $\beta$-Lipschitz continuous for $p_i\in{\cal P}_i$. Accordingly, the minimum singular value of $\nabla g_i(p_i^t)$ is equal to smallest absolute eigenvalue of $[\nabla^2 f_i(g_i(p^t))]^{-1}$, i.e., $\sigma_{\min}(\nabla g_i(p_i^t))=1/L$. This implies that if $\|\hat{\nabla}g_i^t-\nabla g_i(p_i^t)\|\leq 1/(2L)$ holds, then
\begin{align}
    &\sigma_{\min}(\hat{\nabla}g_i^t)=\underset{\|x\|=1}{\min}\|\hat{\nabla}g_i^t x\|\\
    &= \underset{\|x\|=1}{\min}\|\nabla g_i(p_i^t)x+(\hat{\nabla}g_i^t-\nabla g_i(p_i^t))x\|\\
    &\geq \underset{\|x\|=1}{\min}\|\nabla g_i(p_i^t)x\|-\underset{\|x\|=1}{\max}\|(\hat{\nabla}g_i^t-\nabla g_i(p_i^t))x\|\\
    &=1/L-1/(2L)\geq 1/(2L),
\end{align}
which implies that $\|[\hat{\nabla}g_i^t]^{-1}\|=1/\sigma_{\min}(\hat{\nabla}g_i^t)\leq 2L$.
\subsection{Proof of Lemma~\ref{lem:xstepamount}}\label{app:xstepamount}
To bound $\|\hat{x}_i^{t+1}-x_i^t\|$, we will use the following as an auxiliary result:
     \begin{theorem}\cite[Theorem 1.2.1]{schneider2014convex}\label{thm:schneider121}
         Let ${\cal X}$ be a convex and compact set in $\mathbb{R}^d$. Then, the metric projection onto ${\cal X}$ is contracting, that is,
         $$\|\Pi_{\cal X}(x)-\Pi_{\cal X}(y)\|\leq \|x-y\|,~\forall x,y,\in{\mathbb R}^d.$$
     \end{theorem}
     Using the above result, we bound $\|\hat{x}_i^{t+1}-x_i^t\|$ as:
    \begin{align}
        &\|\hat{x}_i^{t+1}-x_i^{t}\| \leq \|\hat{x}^{t+1}-x^{t}\|\\
        &=\|\Pi_{{\cal X}_{\Delta^t}}(x^t+p^t\gamma^t)-\Pi_{{\cal X}_{\Delta^t}}(x^t)+\Pi_{{\cal X}_{\Delta^t}}(x^t)-x^t\|\\
        &\leq \|\Pi_{{\cal X}_{\Delta^t}}(x^t+p^t\gamma^t)-\Pi_{{\cal X}_{\Delta^t}}(x^t)\| {+} \|\Pi_{{\cal X}_{\Delta^t}}(x^t)-x^t\|\\
        &\leq \|p^t\gamma^t\| + \Delta^t\Gamma_{\cal X}\leq M\sqrt{n}\gamma^t+ \Delta^t\Gamma_{\cal X},
    \end{align}
    where we used $\|p_i^t\|=\|\nabla f_i(x_i^t)\|\leq  M$ since $x_i^t\in{\cal X}_i^{\interior}$, and Proposition~\ref{prop:sharp_conv}.

\subsection{Proof of Lemma~\ref{lem:smoothg}}\label{app:smoothg}
The first part of the lemma follows from  the same steps as in Lemma~\ref{lem:gprops} for $p_i\in{\mathbb R}^{d_i}$ instead of $p_i\in{\cal P}_i$, and using $\tilde{f}_i$ and $\tilde{g}_i$ instead of $f_i$ and $g_i$.

Next, we prove the second part of the lemma. For the first statement, given a $p_i\in{\mathbb R}^{d_i}$, suppose that $\tilde{g}_i(p_i)\in{\cal X}_i^\interior$. This implies that there exists $x_i\in{\cal X}_i^{\interior}$ that satisfies $\nabla \tilde{f}_i(x_i)=p_i$. Since $\tilde{f}_i(x_i)=f_i(x_i)$ for $x_i\in{\cal X}_i^\interior$, the same $x_i$ solves the optimization problem in \eqref{eq:priceresponse}, which implies $g_i(p_i) = \tilde{g}_i(p_i)$. Therefore, $g_i(p_i)\in{\cal X}_i^\interior$, which proves $p_i\in{\cal P}_i$ by definition.

To prove the second statement, note that if $p_i\in{\cal P}_i$, then $g_i(p_i)\in{\cal X}_i^\interior$. Since ${\cal X}_i\subseteq \dom f_i$ by Assumption~\ref{ass:feasibleset}, the first order optimality condition of \eqref{eq:priceresponse} implies that there exists $x_i=g_i(p_i)\in{\cal X}_i^\interior$ such that $\nabla f_i(x_i) = p_i$. The same $x_i$ solves the optimization problem \eqref{eq:modifiedresponse}, since $f_i(x_i)=\tilde{f}_i(x_i)$ for $x_i\in{\cal X}_i^\interior$. The optimal solution to \eqref{eq:modifiedresponse} has to be unique due to strong concavity, therefore it must hold true that $\tilde{g}_i(p_i) = g_i(p_i)$.

\makeatletter
\iftoggle{nomsort}{%
    \let\old@@@nomenclature=\@@@nomenclature        
        \newcounter{@nomcount} \setcounter{@nomcount}{0}%
        \renewcommand\the@nomcount{\two@digits{\value{@nomcount}}}
        \def\@@@nomenclature[#1]#2#3{
          \addtocounter{@nomcount}{1}%
        \def\@tempa{#2}\def\@tempb{#3}%
          \protected@write\@nomenclaturefile{}%
          {\string\nomenclatureentry{\the@nomcount\nom@verb\@tempa @[{\nom@verb\@tempa}]%
          \begingroup\nom@verb\@tempb\protect\nomeqref{\theequation}%
          |nompageref}{\thepage}}%
          \endgroup
          \@esphack}%
      }{}
\makeatother
\setlength{\nomlabelwidth}{1.5cm}
\nomenclature{$n$}{Number of users}
\nomenclature{$f_i$}{Utility function of user $i$}
\nomenclature{$\nabla f_i$}{Gradient of $f_i$}
\nomenclature{$f$}{Sum of the $n$ $f_i$'s}
\nomenclature{$\nabla f$}{Gradient of $f$}
\nomenclature{$x_i$}{Resource demand vector of user $i$}
\nomenclature{$x$}{Concatenated resource demand of $n$ $x_i$'s}
\nomenclature{$x^\star$}{Optimal solution}
\nomenclature{$f^\star$}{Optimal objective value}
\nomenclature{$d_i$}{Dimension of $x_i$}
\nomenclature{$\bar{d}$}{Highest dimension of $x_i$ among the users}
\nomenclature{$\dom f_i$}{Domain of $f_i$}
\nomenclature{$\cal X$}{Feasible set}
\nomenclature{${\cal X}_i$}{the set of values that user $i$’s resource demand vector can take in the feasible set ${\cal X}$}
\nomenclature{${\cal X}^\interior, {\cal X}_i^\interior$}{Interiors of sets ${\cal X}$, ${\cal X}_i$}
\nomenclature{$R$}{Upper bound on the diameter of ${\cal X}$}
\nomenclature{$\mu$}{Strong concavity constant for all $f_i$}
\nomenclature{$L$}{Smootness constant for all $f_i$}
\nomenclature{$M$}{Lipschitz constant for all $f_i$}
\nomenclature{$\beta$}{Smoothness constant for all $\nabla f_i$}
\nomenclature{$p_i$}{Resource price vector for user $i$}
\nomenclature{$p$}{Concatenated resource price vector of $n$ $p_i$'s}
\nomenclature{$g_i$}{Price response function of user $i$}
\nomenclature{$g$}{Concatenated price response function of $n$ $g_i$'s}
\nomenclature{$R(T)$}{Regret incurred after $T$ iterations}
\nomenclature{$t$}{Iteration index}
\nomenclature{$p_i^t$}{Price vector of user $i$ at iteration $t$}
\nomenclature{$p^t$}{Concatenated price vector of $n$ $p_i^t$'s}
\nomenclature{$x_i^t$}{Resource demand vector of user $i$ at iteration $t$}
\nomenclature{$x^t$}{Concatenated resource demand vector of $n$ $x_i^t$'s}
\nomenclature{$\bar{\mathcal{B}}(r)$}{Closed Euclidean ball with radius $r$ centered at origin}
\nomenclature{${\cal B}(r)$}{Open Euclidean ball with radius $r$ centered at origin}
\nomenclature{${\cal X}_\Delta$}{Shrunk version of ${\cal X}$ by an amount ${\Delta}$}
\nomenclature{$H_{\cal X}$}{Maximum shrinkage of ${\cal X}$}
\nomenclature{$\hat{x}_i^t$}{Desired resource demand vector of user $i$ at iteration $t$}
\nomenclature{$\hat{x}^t$}{Desired concatenated resource demand vector of $n$ $\hat{x}_i^t$'s}
\nomenclature{$\hat{\nabla}g_i^t$}{Jacobian estimate of user $i$'s price response function at iteration $t$}
\nomenclature{$p_i^{t,s}$}{Resource price for user $i$ at sampling stage of iteration $t$}
\nomenclature{$x_i^{t,s}$}{Resource demand of user $i$ at sampling stage of iteration $t$}
\nomenclature{$\Delta^t$}{Amount of shrinkage of the feasible set ${\cal X}$ at iteration $t$}
\nomenclature{$\gamma^t$}{Step-size of the algorithm at the update stage}
\nomenclature{$\eta^t$}{The amount of price variation at the sampling stage}
\nomenclature{$\tau$}{Constant shift in the denominator of $\gamma^t$ and $\Delta^t$}
\nomenclature{$\mathrm{Sharp}_{\mathcal{X}}$}{Sharpness of ${\cal X}$}
\nomenclature{$\Gamma_{\cal X}$}{Upper bound on the sharpness constant of ${\cal X}$}
\nomenclature{${\cal P}_i$}{Set of prices that induce a resource demand in ${\cal X}_i$ for user $i$}
\immediate\write18{makeindex \jobname.nlo -s nomencl.ist -o \jobname.nls}
\printnomenclature
\newpage
\section*{Nomenclature}

\begin{tabular}{p{3cm}p{8cm}}
$n$ & Number of users \\
$f_i$ & Utility function of user $i$ \\
$\nabla f_i$ & Gradient of $f_i$ \\
$f$ & Sum of the $n$ $f_i$'s \\
$\nabla f$ & Gradient of $f$ \\
$x_i$ & Resource demand vector of user $i$ \\
$x$ & Concatenated resource demand of $n$ $x_i$'s \\
$x^\star$ & Optimal solution \\
$f^\star$ & Optimal objective value \\
$d_i$ & Dimension of $x_i$ \\
$\bar{d}$ & Highest dimension of $x_i$ among the users \\
$\dom f_i$ & Domain of $f_i$ \\
$\cal X$ & Feasible set \\
${\cal X}_i$ & The set of values that user $i$’s resource demand vector can take in the feasible set ${\cal X}$ \\
${\cal X}^\interior, {\cal X}_i^\interior$ & Interiors of sets ${\cal X}$, ${\cal X}_i$ \\
$R$ & Upper bound on the diameter of ${\cal X}$ \\
$\mu$ & Strong concavity constant for all $f_i$ \\
$L$ & Smootness constant for all $f_i$ \\
$M$ & Lipschitz constant for all $f_i$ \\
$\beta$ & Smoothness constant for all $\nabla f_i$ \\
$p_i$ & Resource price vector for user $i$ \\
$p$ & Concatenated resource price vector of $n$ $p_i$'s \\
$g_i$ & Price response function of user $i$ \\
$g$ & Concatenated price response function of $n$ $g_i$'s \\
$R(T)$ & Regret incurred after $T$ iterations \\
$t$ & Iteration index \\
$p_i^t$ & Price vector of user $i$ at iteration $t$ \\
$p^t$ & Concatenated price vector of $n$ $p_i^t$'s \\
$x_i^t$ & Resource demand vector of user $i$ at iteration $t$ \\
$x^t$ & Concatenated resource demand vector of $n$ $x_i^t$'s \\
$\bar{\mathcal{B}}(r)$ & Closed Euclidean ball with radius $r$ centered at origin \\
${\cal B}(r)$ & Open Euclidean ball with radius $r$ centered at origin \\
${\cal X}_\Delta$ & Shrunk version of ${\cal X}$ by an amount ${\Delta}$ \\
$H_{\cal X}$ & Maximum shrinkage of ${\cal X}$ \\
$\hat{x}_i^t$ & Desired resource demand vector of user $i$ at iteration $t$ \\
$\hat{x}^t$ & Desired concatenated resource demand vector of $n$ $\hat{x}_i^t$'s \\
$\hat{\nabla}g_i^t$ & Jacobian estimate of user $i$'s price response function at iteration $t$ \\
$p_i^{t,s}$ & Resource price for user $i$ at sampling stage of iteration $t$ \\
$x_i^{t,s}$ & Resource demand of user $i$ at sampling stage of iteration $t$ \\
$\Delta^t$ & Amount of shrinkage of the feasible set ${\cal X}$ at iteration $t$ \\
$\gamma^t$ & Step-size of the algorithm at the update stage \\
$\eta^t$ & The amount of price variation at the sampling stage \\
$\tau$ & Constant shift in the denominator of $\gamma^t$ and $\Delta^t$ \\
$\mathrm{Sharp}_{\mathcal{X}}$ & Sharpness of ${\cal X}$ \\
$\Gamma_{\cal X}$ & Upper bound on the sharpness constant of ${\cal X}$ \\
${\cal P}_i$ & Set of prices that induce a resource demand in ${\cal X}_i$ for user $i$ \\
\end{tabular}
\end{document}